\documentclass[12]{article}
\usepackage[left=2.8cm,top=2.8cm,right=2.8cm,bottom=2.8cm]{geometry}
\usepackage{cite}
\usepackage{epsfig}
\usepackage{enumerate}
\usepackage{amssymb}
\usepackage{tikz}
\usepackage{pgfplots}
\usetikzlibrary{matrix,arrows,decorations.pathmorphing}
\usepackage{amsmath}
\usepackage{amsthm}
\usepackage{color}
\usepackage{parskip}
\usepackage{indentfirst}
\newcommand{\octant}{orthant}
\newcommand{\octants}{orthants}
\newcommand{\arctanh}{\text{arctanh}}
\newcommand{\sech}{\text{sech}}
\newcommand{\diagonal}{\text{diag}}
\newcommand{\sign}{\text{sign}}
\newtheorem{theorem}{Theorem}[section]

\newtheorem{cor}{Corollary}[section]
\newtheorem{prop}{Proposition}[section]
\newtheorem{lem}{Lemma}[section]
\newtheorem{rem}{Remark}[section]

\theoremstyle{definition}
\date{}

\title{Introduction}
\begin{document}

\begin{center}

\textbf{ Stability of Discrete Time Recurrent Neural Networks and Nonlinear optimization problems} \end{center}
 \begin{center} Dr. Nikita Barabanov, and Jayant Singh\end{center}

\textbf{Abstract} We consider the method of Reduction of Dissipativity Domain to prove global Lyapunov stability of Discrete Time Recurrent Neural Networks. The standard and advanced criteria for Absolute Stability of these essentially nonlinear systems produce rather weak results. The method mentioned above is proved to be more powerful. It involves a multi-step procedure with maximization of special nonconvex functions over polytopes on every step. We derive conditions which guarantee an existence of at most one   point of local maximum for such functions over every hyperplane. This nontrivial result is valid for wide range of neuron transfer functions.  
\section{Introduction and Problem Setting}
 
\noindent In this paper, we study existence of points of local maxima  for function $f(x)=\sum_{i=1}^{n}c_i \phi(x_i),$ where $\phi(\cdot)$ is a nonlinear function, over a hyperplane. This problem arises in stability analysis of nonlinear dynamical systems \cite{Nikita}, for example Recurrent Neural Networks (RNN). A typical RNN can be described by the following systems of equations;

\begin{align}
\label{eq1}
x_{1}^{k+1} &= \phi(W_1 x_{1}^{k}+V_n x_{n}^{k}+b_1),\nonumber \\
x_{2}^{k+1} &= \phi(W_2 x_{2}^{k}+V_{1} x_{1}^{k+1} +b_2),\nonumber \\
\cdots \nonumber \\
x_{n}^{k+1} &= \phi(W_n x_{n}^{k}+V_{n-1}x_{n-1}^{k+1}+b_n),
\end{align}

\noindent where  $n$ is the number of layers, $\phi(\cdot)$ is the activation function,
$x_{j}^{k}$ is the state vector of the layer $j$ at time step $k$,
$W_j $ and $V_j$ are fixed weight matrices,  and
$b_j $ is a fixed vector representing bias. 

\indent Stability of RNN has been addressed extensively in literature. In \cite{Sukyens}, a stability criterion has been developed using $NL_q$ approach. A typical $NL_q$ system( without external inputs), is of the form,
\begin{equation}
\label{eqyyyy}
p_{k+1}=P_1 Q_1 P_2 Q_2 \ldots P_q Q_q p_k
\end{equation}
where $p_k \in \mathbb{R}^n$, $P_i=(\diagonal{\overline{p_j}})_{j=1}^{n}$,  and  $Q_i$is a constant matrix. Here $\overline{p_j}$ depends on $p_k$ continuously. The problem under consideration is to check stability of system \eqref{eqyyyy}, with matrices $P_i$ satisfying the relation $\Vert P_i\Vert \leq 1.$  The stability criterion using $NL_q$ approach says that if there exists diagonal positive definite matrices $D_j$ such that $\Vert D_j Q_j D_{j+1}^{-1}\Vert <1 $ for all $j=1,\ldots ,q(mod \;q)$, then the system \eqref{eqyyyy} is globally asymptotically stable. Using a suitable method \cite{NBarabanov}, the RNN defined in \eqref{eq1} can be transformed to form \eqref{eqyyyy}. Therefore the above criterion can be used to check stability of systems of form \eqref{eq1}. The $NL_q$ approach gives sufficient conditions for stability of nonlinear systems. However, there exist  nonlinear stable systems, for which the $NL_q$ stability criterion is not satisfied. These nonlinear systems, for example, RNN, have shown promise in various applications \cite{Feldkamp}.  
 
\indent  Another stability criterion was developed using theory of absolute stability (\cite {NBarabanov}, \cite{NikitaBarabanov}, \cite{Yakub}, \cite{Narendar}, \cite{Aiz}).  A system to be analyzed for stability using this approach  should be written in the automatic control form:
\begin{equation} 
\begin{array}{lcr}
x^{k+1} &=& A x^{k}+B \psi^{k},\\
\sigma^{k}&=&C x^{k},\\
\psi_{i} &=& \phi _{i}(\sigma_{i}),i=1 \cdots m,
\end{array}
\label{eq2}
\end{equation} 
where, $A,B,C$ are matrices of suitable size,
$\psi^{k}= ( \psi_1,  \ldots  ,\psi_m )$ is input vector at step $k$,
$\sigma^{k}=( \sigma_1,  \ldots , \sigma_m )$ is the output vector at step $k$, and 
$\{\phi_{i}(\cdot)\}_{i=1}^{m}$ are nonlinear functions.

\indent Before analyzing the stability of system \eqref{eq1} using theory of Absolute stability, it needs to be transformed to \eqref{eq2}. 
State Space  Extension method has been introduced \cite{NBarabanov} to transform RNN to  \eqref{eq2}.

\indent One of the  significant contribution of theory of absolute stability is the frequency domain criterion (\cite{Kalman},\cite{Yakubo},\cite{Yak}). Frequency domain criterion gives necessary and sufficient condition for existence of  quadratic Lyapunov function for class of systems \eqref{eq2} with functions $\phi(\cdot)$ satisfying given local quadratic constraint. One of the most common constraints used for stability analysis of nonlinear systems  is sector constraint, and the corresponding stability criterion is known as   circle criterion. It has been shown in \cite{NBarabanov} that stability criterion given by $NL_q$ approach is weaker than the circle criterion.

\indent The circle criterion gives sufficient condition for stability of nonlinear systems, with nonlinear function $\phi(\cdot)$ satisfying sector constraint. It only utilizes the fact that the nonlinear function $\phi(\cdot)$ satisfies a given sector condition. It might happen that given a sector, defined by function $\phi(\cdot),$  there exists a nonlinear function satisfying sector condition, such that the  corresponding system is unstable. Additional information about the nonlinear function can be used to check stability of nonlinear systems of particular kind, for example RNN. A modified stability criterion using additional information about the nonlinear function,( e.g. monotonicity) has been developed in \cite{NBarabanov}. But this criterion has been shown to be essentially sufficient for systems with large number of nonlinear functions.  In addition, this criterion is not applicable to  some practically stable systems, for instance RNN.

\indent The stability criterion given by theory of absolute stability ( \cite{NBarabanov},\cite{Yakubovich}, \cite{Molinari},\cite{Lankaster}) checks necessary and sufficient conditions for existence of Lyapunov functions of a particular kind (e.g. quadratic forms). But there exists stable systems, for which quadratic Lyapunov functions do not exist. An alternative stability criterion has been proposed in \cite{Nikita}. 

Consider the system 
\begin{equation}
\label{eqyyy}
x_{k+1}=\phi(x_k)
\end{equation}

 \indent Let $D_0$ denote the whole space of vector $x_k.$ Suppose there exists sets $\{D_k\}$ such that $D_{k+1} \subset D_k, \phi(D_k)\subseteq D_{k+1}$.  If $\{D_k\} \rightarrow 0$ (in Hausdorff metric), as $k \rightarrow \infty$, then obviously system \eqref{eqyyy} is globally asymptotically stable.  This approach is known as reduction of dissipativity domain. \\
 \indent  In order to implement this approach, the sets $D_k$ need to be defined.   A possible choice of $D_k$ is given by 
 \begin{equation*}
D_{k+1}=\{x \in D_k: f_{k+1,j}(x)\leq \alpha_{k+1,j},j=1\ldots m_{k+1}\}
\end{equation*}
\indent where $m_k$ is the number of constraints at step $k,
 f_{k,j}$ is a  function,  and 
 $\alpha_{k+1,j}=\max_{x \in D_k}f_{k,j}(\phi(x)).$ \\  The set $D_{k}$ is characterized by the set of pairs $(f_{k,j},\alpha_{k,j})$ where $j \in \{1 \ldots m\}$. A possible choice of $f_{k,j}(\cdot)$ is linear functions.   Then $D_k$ takes the shape of a polytope. It has been shown (reference) that if system \eqref{eqyyy} has a convex Lyapunov function, then there exists linear functions $f_{k,j}$ such that $\{D_k\}\rightarrow 0$.  
 
 \indent  The set $D_k$ is constructed by computing the value  $\alpha_{k,j}$ for every $j.$ Since the function $\phi(\cdot)$ is nonconcave over the set $D_k$, it can have multiple points of local maxima.  At every step $k,$ the points of local maxima for the function $f(\phi(\cdot))$ need to be computed.\\
 \indent Consider a single layer RNN with zero bias. Using substitution $y=Wx,$ it can be expressed as 
 
 \begin{equation}
 \label{eqzzz}
 y_{k+1}=W\phi(y_k)
\end{equation}

 \indent For the case of RNN in \eqref{eqzzz},the function $f(\phi(\cdot))$ is given  by the inner product $f(x):=\left\langle l_j,W\phi(x)\right\rangle $. We need to find points of local maxima for $f(\cdot)$  over polytopes defined by matrix of constraints, $L=col(l_1, l_2,\ldots,l_m).$ It has been seen that, in all the  cases, the  function $ f(\cdot)$ has points of local maxima on the boundary  of the polytope.  We will first locate the points of local maxima for $f(\cdot)$ on an arbitrary hyperplane. The subject of this paper is the solution to the following problem.\\

\indent  \underline{Problem Setting}: Consider the hyperplane, $P=\{x:l^{T} x=b\} $ where $l$ is a unit normal vector and $b \in \mathbb{R}.$ How many points of local maxima does the  function $f(x)=\sum_{i=1}^{n}c_i \phi(x_i),c_i \neq 0$ for all $i,$ have on $P?$ Here $\phi(\cdot)$ is standard neuron transfer function. 

\indent This paper is organized as follows. In section 2, we will develop necessary and sufficient conditions for existence of points of local maxima. In section 3, some assumptions regarding function $\phi(\cdot)$ will be listed.   Section 4 gives the possible location of points of local maxima. Then we will talk about number of points of local maxima in main orthant and side orthant.  We will conclude with the main result of this paper.

\section{Identify the points of local maxima}

\indent We will develop the necessary and sufficient condition for a critical point to be a point of local maximum for function $f(\cdot)$ over hyperplane $P.$ Then $D=\frac{\partial ^ {2}f}{\partial x^2}\mid_{x=x_0}=\diagonal(d_j)_{j=1}^{n}$ is Hessian matrix, where $d_j=c_j \phi''(x_j).$ Let $K:=(I-\frac{ll^{T}}{\Vert l \Vert ^{2}})D(I-\frac{ll^{T}}{\Vert l \Vert ^2})$  denote the  projection matrix.
\begin{theorem}
\label{thm1}
 Suppose, $x_0$ is a critical point of $f(\cdot)$ over $P$ (i.e. $l$ is parallel to $\overrightarrow{\triangledown f}(x_0)$).  
Then $x_0$ is a point of local maximum of $f(\cdot)$ over P  only if $K \leq 0.$ Moreover, if $K$ has $n-1$ negative eigenvalues and one zero eigenvalue then $x_0$ is a point of local maximum.
\begin{proof}
\textit{Necessity}: Consider the Taylor expansion  for $f(\cdot)$ in some neighborhood of $x_0.$\\
$$f(x)=f(x_0)+\left\langle \frac{\partial f}{\partial x}\vert_{x=x_0}, x-x_0\right\rangle+\frac{1}{2}\left\langle x-x_0, \frac{\partial^{2}f}{\partial x^2}\vert_{x=x_0} (x-x_0)\right\rangle+o(\Vert x-x_0\Vert ^2)$$
Since $x\in P,$ we have  $l^{T} x=b$ and $l^{T} x_0=b,$ hence $x-x_0$ is orthogonal to $ l.$ Using the fact that, $\overrightarrow{\triangledown f}(x_0)$ is parallel to $ l,$ we get  $\overrightarrow{\triangledown f}(x_0)^T (x-x_0)=0.$ 

\indent Moreover since $x_0$ is a point of local maximum, we obtain $\left\langle x-x_0, \frac{\partial^{2}f}{\partial x^2}\vert_{x=x_0} (x-x_0)\right\rangle \leq 0$  for all $x$ such that  $l^T(x-x_0)=0$. Therefore 
\begin{equation}
\label{eq3}
  y^{T} D y \leq 0 
  \end{equation}
   for all $ y$ such that $ l^{T}y=0.$

 \indent Pick $z \in \mathbb{R}^n.$  Define   $y:=(I-\frac{l l^{T}}{\Vert l \Vert ^{2}})z$. Then $y \in P.$ Next we will show that $K\leq 0.$

\indent Pick $z$ arbitrary. Then
\begin{align}
z^{T}K z &= z^{T}\Big(I-\frac{ll^{T}}{\Vert l \Vert ^{2}}\Big)D\Big(I-\frac{ll^{T}}{\Vert l \Vert ^2}\Big)z \nonumber\\
&= z^{T}\Big(I-\frac{ll^{T}}{\Vert l \Vert ^{2}}\Big)^{T}D \Big(I-\frac{ll^{T}}{\Vert l \Vert ^{2}}\Big)z\nonumber.
\end{align}
\indent Hence, $z^{T}Kz=\Big((I-\frac{ll^{T}}{\Vert l \Vert ^{2}})z\Big)^{T}D \Big((I-\frac{ll^{T}}{\Vert l \Vert ^{2}}) z\Big)= y^{T} D y$. Using equation$\eqref{eq3},$ we have  $y^{T}Dy\leq 0.$ Therefore, $K\leq 0.$

\indent \textit{Sufficiency}: Suppose $K$ has $(n-1)$ negative eigenvalues and one zero eigenvalue. This implies that $z^ T K z \leq 0$ for all $z \in \mathbb{R}^n.$ We will show that if $z \in  \{x:l^{T}x=0\},$ then $z^T K z<0.$ 

\indent There exists orthonormal basis $\{v_1,v_2,\ldots,v_{n-1}\}$ of P consisting of eigenvectors of matrix K, with eigenvalues $\{\lambda_1,\lambda_2,\ldots,\lambda_{n-1}\}$. Since $v_0=l$ is the eigenvector with zero eigenvalue, we get $\lambda_i <0$ for all $i \in \{1\ldots n-1\}.$

\indent Assume $z \in \{x:l^{T}x=0\}.$ Then $z=\sum_{j=1}^{n-1}p_j v_j, p_j \in \mathbb{R}$ and,
\begin{align}
z^TKz &= \left\langle \sum_{i=1}^{n-1}p_i v_i,\sum_{i=1}^{n-1}p_i Kv_i\right\rangle\nonumber
=\left\langle \sum_{i=1}^{n-1}p_i v_i,\sum_{i=1}^{n-1}p_i \lambda_i v_i\right\rangle \nonumber\\
&=\sum_{i=1}^{n-1}\lambda_i p_{i}^{2}\nonumber
<max(\lambda_i)\sum_{i=1}^{n-1}p_{i}^{2}
<0.\nonumber
\end{align}

\indent Pick $x$ in neighborhood of $x_0$ on P. Then $(x_0 -x) \perp l.$ Put $z=x_0 -x.$ We obtain $(x_0 -x)^{T}K (x_0- x)<0.$ Using definition of matrix $K,$ we get $(x_0-x)^{T}K(x_0 -x)=(x_0 -x)^{T}D(x_0 -x)<0.$ Since $x_0$ is critical point, $\left\langle \frac{\partial f}{\partial x}\vert_{x=x_0}, x-x_0\right\rangle=0.$ Using Taylor expansion for $f(\cdot),$ we obtain $f(x)<f(x_0)$ in some neighborhood of $x_0$ on P. Therefore $x_0$ is a point of local maximum for $f(\cdot)$ on hyperplane P.
\end{proof}
\end{theorem}

\indent In the following section, we will list some assumptions about the function $\phi(\cdot).$ These assumptions  will be used to show the main result  of this paper.

\section{Assumptions about  Cost function}

\underline{Notation}: The following notation will be followed, unless specified.

\indent $\psi'(s):=\frac{d\psi}{dx}\vert_{x=s},\psi_{\beta }(\beta q):=\frac{d}{d\beta}(\psi(\beta q)), h_{\beta}(\beta,q_j,q_n):=\frac{\partial}{\partial \beta} (h(\beta,q_j,q_n))$,where $\psi(\cdot), $ and $h(\cdot)$ are functions which will be defined later.

\indent \textit{Assumption 1}: $\phi(\cdot) \in C^2,\phi(-x)=-\phi(x),\phi'(x)>0, x\phi''(x)<0,$ for all $x\neq 0,$ and $\lim_{x \rightarrow \infty}\phi(x)<\infty.$ 
$x \phi''(x)<0$ implies that $\phi'(x)$ is decreasing function for all $x>0.$ Hence $\phi'(\cdot)^{-1}$ exists. Denote $\phi'(\cdot)^{-1}=\psi(\cdot).$ We get $\psi:(0,\phi'(0)]\rightarrow [0,\infty). $ In addition, $\psi'(x)<0$ for all $x.$ 

\indent \textit{Assumption 2}: $x(\ln \vert\psi'(x)\vert)'$  is a monotonically increasing function of $x.$
This implies that   $ \frac{d}{d\beta}\Big(\frac{\psi'(\beta p)}{\psi '(\beta q)}\Big)>0$, where $p>q$.

\indent \textit{Assumption 3}: Set $h(\beta, q_j,q_n) =\frac{\psi'(\beta q_j)}{\psi'(\beta q_n)}.$ Then $\frac{\partial}{\partial \beta}\Big[ \frac{h_{\beta}(\beta, q_j,q_n)}{h_{\beta}(\beta, q_l,q_n)}\Big]$ is sign definite, where $ q_j <q_n <q_l.$

 \indent \textit{Assumption 4} : For all $p>q,$ we have $\frac{d}{d\beta}\Big(\frac{\psi(\beta p)}{\psi(\beta q)}\Big)<0.$
 
 \indent \textit{Assumption 5}: For all $x>0, $ we have $\frac{d}{dx}\Big(x \frac{d}{dx}\Big(\frac{\psi(x)}{x\psi'(x)}\Big)\Big)\geq 0.$\\

 \section{Possible Locations of Points of Local Maxima}

\indent In this section, we will use above assumptions to locate the possible locations of points of local maxima for function  $f(\cdot)$ over hyperplane $P.$ Recall $f(x)=\sum_{i=1}^{n}c_i \phi(x_i)$ and   $P=\{x:l^{T}x=b\}$, where $l=(l_j)_{j=1}^{n}$ is the normal vector.

\indent First, we change basis in order to get $c_j >0$ for all $j \in \{1 \ldots n\}.$ Suppose $c_{j}=0$ for some $j \in \{1\ldots n\}$, then the corresponding term in the sum is zero, and we obtain $f(x)=\sum_{k=1}^{n-1}c_k x_k.$ The problem is reduced to similar problem of dimension $n-1.$ Without loss of generality,  we assume that  $c_j \neq 0$ for all $j \in \{1\ldots n\}.$ Next assume that $c_{j0}<0 $ for some $j0 \in \{1\ldots n\}.$ Using assumption (1), $\phi(\cdot)$ is odd function. This implies that $c_{j0} \phi(x_{j0})=-c_{j0}\phi(-x_{j0}).$  Hence, if $c_{j0} <0 $ then  replacing $c_{j0}$ by $-c_{j0}$, and $e_{j0}$ by $-e_{j0}$, the function $f(x)$ remains unchanged, and coefficient $c_{j0}>0.$ Without loss of generality, we assume $c_j >0 $ for all $j \in \{1 \ldots n\}.$

\indent Next, we consider the signs of the components of vector $l.$ If $l_{j0}=0 $ for some $j0 \in \{1 \ldots n \}$, then $x_{j0}$ can be increased arbitrarily, still $\sum_{j=1}^{n}l_j x_j$ remains unchanged. Hence, the function $f(x)$ does not have  a point of local maximum on $P.$ Therefore we can assume that $l_j \neq 0$ for all $j \in \{1 \ldots n\}.$  Next, assume  that  there exists $j0, j1 \in \{1\ldots n\}$ where $j0 \neq j1$, such that $l_{j0}<0<l_{j1}.$  Then we will increase $x_{j0}$ and  $x_{j1}$ such that $\sum_{j=1}^{n}l_j x_j$ is unchanged. It is easy to see that the function $f(x)$ is increasing on $P.$ The function $f(x)$ does not have  a point of local maximum on $P.$ Hence, we can assume that for any $j0,j1 \in \{1 \ldots n \},$ where $j0 \neq j1$ the product  $l_{j0}l_{j1}$ is positive. Suppose $l_j<0$ for all $j \in \{1\ldots n\}.$ Then  we replace $b $ by $-b$ and vector $l$ by $-l.$ Without loss of generality, we  assume that $l_j>0 $ for all $j \in \{1\ldots n\}.$

\begin{theorem}
\label{thm2}
Suppose $x_0$ is a critical point for function $f(\cdot).$ Then, $x_0$ is a point of local maximum only if  the \octant \;  has at most one negative coordinate (i.e. $x_j <0$ for at  most one $j$, where $j \in \{1 \ldots n\}$).
\indent \begin{proof} Notice that $D=\diagonal(d_j)_{j=1}^{n}$, and $l$ denotes the normal vector. The characteristic polynomial of $K$ is given by $\det(\lambda I-K)=0$, where $\lambda$ denotes eigenvalue of matrix $K.$ First we will compute $\det(\lambda I-K).$ 
\begin{align}
\det(\lambda I-K)&= \det\Big (\lambda I-(I-\frac{ll^{T}}{\Vert l \Vert ^{2}})D(I-\frac{ll^{T}}{\Vert l \Vert ^2})\Big)\nonumber\\
                &= \det \Big(\lambda I- D(I-\frac{ll^{T}}{\Vert l \Vert ^{2}})(I-\frac{ll^{T}}{\Vert l \Vert ^{2}})\Big)\nonumber\\
                 &= \det\Big(\lambda I-D\Big)\det\Big(I+\frac{l^{T}(\lambda I-D)^{-1}Dl}{\Vert l\Vert ^{2}}\Big ).\nonumber
\end{align}

\indent using Sylvester identity.

\indent Hence,
\begin{align}
\det(\lambda I-K) &= \det\Big(\lambda I-D\Big)\Big(1+ \frac{l^{T}(\lambda I-D)^{-1}Dl}{\Vert l\Vert ^{2}}\Big)\nonumber\\
                   &= \Big(\Pi_{i=1}^{n}(\lambda- d_i)\Big)\cdot \sum_{i=1}^{n}\frac{\lambda l_{i}^{2}}{\Vert l \Vert ^{2}(\lambda-d_i)}.
\label{eq4}
\end{align}

It follows from equation (\ref{eq4}) that $\lambda=0$ is an eigenvalue of $K.$ Denote $g(\lambda)=\sum_{i=1}^{n}\frac{\lambda l_{i}^{2}}{\Vert l \Vert ^{2}(\lambda-d_i)}.$

\indent  We can see that $g(\lambda)$ has vertical asymptotes at $d_j$. Suppose all $d_{j}$'s are distinct. Then, we can arrange them as $d_1 <d_2<\ldots<d_n.$  It is easy to see that  for all $d_j>0,$ \; $ \lim_{\lambda \rightarrow d_{j}^{+}}g(\lambda)=\infty$ and $\lim_{\lambda \rightarrow d_{j}^{-}}g(\lambda)=-\infty$. If $d_j<0,$ then $\lim_{\lambda \rightarrow d_{j}^{+}}g(\lambda)=-\infty$ and $\lim_{\lambda \rightarrow d_{j}^{-}}g(\lambda)=\infty$. Since function $g(\cdot)$ is continuous in $(d_{j},d_{j+1})$, there exists a root, $\lambda_{j}$ of function $g(\cdot)$ in this open interval. Number $\lambda_{j}$ is an eigenvalue of $K$ for $j=1,\ldots,n-1$. Thus, 
$\{0,\lambda_{1},\ldots,\lambda_{n-1}\}$ are all eigenvalues of matrix $K$.
\indent Now consider the general case. We order the values of $d_{j}$: $d_{1}\leq d_{2}\ldots,d_{n}$. If $d_{j}<d_{j+1}$ then the proof is same as  above.  If $d_{j-1}\neq d_{j}=d_{j+1}=\ldots d_{j+k}\ne d_{j+k+1}$, then there are $k$ eigenvalues $\lambda_{j}=\lambda_{j+1}=\ldots=\lambda_{j+k-1}$ of $K$ at point $d_{j}$. 
Together with the zero eigenvalue the set of such numbers $\{\lambda_{j}\}$, $j \in \{1,2,\ldots,n-1\}$ is the set of all eigenvalues of $K$.

\indent \underline{Claim 1}: $d_j >0 $ for at most one $j$ where $j \in \{1 \cdots n\}.$\\
\indent  Suppose there exists $ d_j, d_k,$where $j \neq k$  such that  $d_k>d_j>0.$ Then $g(\lambda)=0 $ for some $\lambda \in (d_j,d_k).$ This implies that matrix $K$ has a positive eigenvalue, but $K\leq 0$ (Theorem \ref{thm1}). Hence, we arrive at a contradiction. Claim 1 proved.

\indent \textit{Case 1}: Suppose $d_j \leq   0$  for all $j.$  Using definition of $d_j,$ and assumption 1, we obtain  $x_j \geq 0$ for all $j.$  

\indent \textit{Case 2}: Suppose $d_j >0 $ for some $j.$  It can be analyzed in a similar manner to Case 1. We obtain that $x_j<0$ for some $j.$\\
\indent Combining the results of Case 1 and Case 2,  we obtain that a  stationary point, $x_0$, is a point of local maximum only if  the \octant \; has at most one negative coordinate.
\end{proof}
\end{theorem}
\indent Using Theorem \ref{thm2}, we can deduce that function $f(x)$ can have points of local maxima in main \octant \; or side \octant \; with at most one coordinate negative. Next, we will  show that $f(x)$ has at most one point of local maximum in main \octant.

\section{Points of Local Maxima in Main Orthant}

\indent In this section, we will show that\; $f(x)$ has at most one point of local maximum in main \octant \;(i.e. $x_j>0$ for all $j \in \{1\ldots n\}$).
\begin{prop}
\label{prop5}
The function $f(x)= \sum_{i=1}^{n} c_i \phi(x_i)$ has at most one  point of local maximum in the main \octant.
\begin{proof}
This is obvious, since $f(\cdot)$ is concave over the main orthant.
\end{proof}
\end{prop}

\indent Next, we will present the necessary and sufficient condition for existence of point of local maximum in main orthant.

\begin{prop}
\label{cor1}
A point $x_0$ is  a point of local maximum in main \octant \; if and only if there exists $\beta \in \mathbb{R}^{+}$such that $\sum_{j=1}^{n}l_j \psi(\beta \frac{l_j}{c_j})=b$, where $\beta\frac{l_j}{c_j}<\phi'(0)$ for all $j$, and  $b$ is such that $l'x_{0}=b.$
\begin{proof} \textit{Necessity}:
 Suppose that  $x_0$ is a point of local maxima in the main \octant \; and $x_0 \in P,$ the hyperplane. Then $x_0$ is a critical point.  We obtain $\overrightarrow{\triangledown {f}}(x_0)=\beta l$ for some $\beta \in \mathbb{R}.$  Since $l_j >0 $ for all $j,$ and $\frac{\partial f}{\partial x_{j}}=c_j \phi'(x_j)>0,$ we obtain that $\beta \in \mathbb{R^{+}}.$

\indent Hence, $\phi'(x_{0}^{j})=\beta \frac{l_j}{c_j}$ for all $ j \in \{1 \ldots n\},$ which in turn implies that $ x_{0}^{j}=\psi(\beta q_j)$ where $q_j:=\frac{l_j}{c_j},$ and $\psi:=(\phi'(\cdot))^{-1}.$ Since $x_0$ lies on hyperplane $P,$ we obtain $b=l'x_0=\sum_{j=1}^{n}l_j \psi(\beta q_j).$
  for some $\beta \in \mathbb{R}^{+}.$

\indent \textit{Sufficiency}: Assume that  there exists $ \beta \in \mathbb{R}^+$ such that $b=\sum_{j=1}^{n}l_{j}\psi(\beta \frac{l_j}{c_j}).$ Using definition of $\psi(\cdot),$ we get $b=\sum_{j=1}^{n}l_{j}(\phi')^{-1}(\beta \frac{l_j}{c_j}).$ Denote  $x_{0}^{j}:=(\phi')^{-1}(\beta \frac{l_j}{c_j}).$  This implies
$ \phi '(x_{0}^{j})=\beta \frac{l_j}{c_j}$
 for all $j \in \{1 \ldots n\}.$  Therefore, $ \overrightarrow{\triangledown }f(x_0)=\beta l.$ Hence $x_0$ is a stationary point in main \octant.  Since $f(x)$ is concave over the main orthant, $x_0$ is a point of local maximum in the main \octant.
\end{proof}
\end{prop}

\indent In Theorem \ref{thm2} we saw that a critical point $x_0 $  for function $f(x)=\sum_{i=1}^{n} c_i \phi(x_i)$ can be a point of local maximum only if it is lying in  \octant \;  with at most one negative component. In this section, we developed the necessary and sufficient conditions for existence of local maximum in main \octant. In addition, we also showed that $f(x)$ can have at most one point of local maximum in main \octant. Next we will show similar result for the case of side \octant \; with one negative component. In later sections we will use \textit{side orthant} to denote side orthant with one negative component.  

\section{ Points of  local maxima in side orthant}
\indent In this section, we will show that the function $f(x)$, defined on the hyperplane, has at most one point of local maximum in side \octant. Here we have shown the result for the side \octant, with last component negative. Other cases can be analyzed similarly.\\
\indent First, we will develop the necessary and sufficient conditions for existence of local maximum in side \octant \;(i.e $x_j>0$ for all $j\in \{1  \ldots n-1\}, x_n<0$). Suppose $x_0$ is a critical point for function $f(\cdot)$ in the side \octant \; (i.e $x_{0}^{j}>0$ for all $j\in \{1  \ldots n-1\}, x_{0}^{n}<0$). Then, there  exists $\beta \in \mathbb{R}^{+}$ such that $c_j \phi'(x_{0}^{j})=\beta l_j$ for all $j \in \{1\ldots n\}$. This implies $x_{0}^{j}=\sign(x_{0}^{j})(\phi ')^{-1}(\beta \frac{l_j}{c_j})=\sign(x_{0}^{j})(\phi ')^{-1}(\beta q_j) $ for all $j \in \{1 \ldots n\}.$ Recall that $d_j=c_j \phi''(x_j)$. Using  assumption 1, we obtain  $d_1 \leq d_2 \leq \cdots \leq d_{n-1}<0<d_n.$ \\
\indent Denote $g_{1}(\beta):=l'x_0=\sum_{j=1}^{n-1}l_j (\phi ')^{-1}(\beta q_j)-l_n (\phi ')^{-1}(\beta q_n)$, where $q_j =\frac{l_j}{c_j}.$ Using the notation $(\phi ')^{-1}:=\psi,$ we obtain

\begin{equation}
g_{1}(\beta)=\sum_{j=1}^{n-1}l_j \psi(\beta q_j)-l_n \psi(\beta q_n)
\label{eq5}
\end{equation}

\begin{theorem}
\label{thm6}
Let $g(\lambda)=\sum_{i=1}^{n}\frac{\lambda l_{i}^{2}}{\Vert l \Vert ^{2}(\lambda-d_i)}.$
If $g'(0)<0$ and $g_{1}(\beta) =b,$ for some $\beta \in \mathbb{R}^{+}$, then $x_0$ is a point of local maximum for the function $f(\cdot)$, over the hyperplane P. Moreover, $x_0$ is a point of local maximum only if $g'(0)\leq 0 $ and $g_{1}(\beta) =b,$ for some $\beta  \in \mathbb{R}^{+}$.
\begin{proof}
\indent \textit{Sufficiency}: It is easy to see that $\lim_{\lambda \rightarrow d_{n}^{-}} g(\lambda)=-\infty$ and $\lim_{\lambda \rightarrow d_{n-1}^{+}}g(\lambda)=-\infty.$ We saw earlier that there exists unique $ \lambda \in (d_{j},d_{j+1})$ for all $j \leq n-2$ and $\lambda $ is eigenvalue of $K.$ Thus matrix $K$ has  $n-2$ negative eigenvalues. \\
\indent So we have 2 roots of $g(\lambda) $ in the interval $(d_{n-1},d_{n}).$  One of the roots  is 0. Denote the other root as $\hat{\lambda}.$ Notice that $d_{n-1}<0 <d_n.$
 This implies that $\hat{\lambda}$ can be negative or positive. But under the assumption that $g'(0)<0, $ we have $\hat{\lambda}<0.$ Moreover, $\hat{\lambda} \in (d_{n-1},0).$ Hence we conclude that matrix $K$ has $n-1$ negative eigenvalues and a zero eigenvalue. Using theorem \ref{thm1},  $x_0$ is a point of local maximum.

\indent \textit{Necessity}: Suppose that $x_0$ is a point of local maximum. Then $K \leq 0.$ This implies $g'(0)\leq 0 $ and 
$b= \sum_{j=1}^{n-1}l_j \psi(\beta q_j)-l_n \psi(\beta q_n)$ for some $\beta \in \mathbb{R}^{+}$.
\end{proof}
\end{theorem}
\indent In theorem \ref{thm6} we developed the necessary and sufficient conditions for existence of point of local maximum in side \octant. These conditions can be rewritten in terms of function $g_{1}(\beta).$ To this end, we need the following equivalent statement.

\begin{lem}
\label{lem4}
Consider function $g(\lambda)$ as defined above. Then $g'(0)<0 $ if and only if  $g_{1}'(\beta)>0.$

\begin{proof}\textit{Necessity}: We can rewrite $ g'(0)$ as 
$$g'(0)=\sum_{j=1}^{n}-\frac{l_{j}^{2}}{d_{j} \Vert l\Vert ^{2}}=\frac{1}{\Vert l\Vert^2}\sum_{j=1}^{n}\frac{l_{j}^{2}}{c_j \phi''(\sign(x_{0}^{j})\psi(\beta \frac{l_j}{c_j}))},$$

\indent where, $d_j=c_j \phi''(x_{0}^{j})=c_j \phi''(\sign(x_{0}^{j})\psi(\beta \frac{l_j}{c_j}))$, and $x_0$ is critical point for $f(\cdot).$ Therefore,
\begin{align} 
g'(0)&=-\frac{1}{\Vert l \Vert ^2}\Big(\sum_{j=1}^{n-1}\frac{l_{j}q_j}{ \phi''(\psi(\beta q_j))}+\frac{l_{n} q_n}{\phi''(-\psi(\beta q_n))}\Big)\nonumber\\
&=\frac{1}{\Vert l\Vert^2}\Big(\frac{l_{n}q_n}{ \phi''(\psi(\beta q_n))}-\sum_{j=1}^{n-1}\frac{l_{j}q_j}{ \phi''(\psi(\beta q_j))}\Big).\nonumber
\end{align}

since $\phi''(\cdot)$ is odd.

\indent Hence $g'(0)<0$ implies 
\begin{equation}
\label{eq6}
\frac{q_n l_n}{ \phi''(\psi(\beta q_n))}<\sum_{j=1}^{n-1}\frac{q_j l_j}{ \phi''(\psi(\beta q_j))}
\end{equation}

\indent Next, we want to express $g_{1}'(\beta)$ in a form similar to $g'(0).$ To this end, we need an auxiliary result.

\indent \underline{Claim 2}: $\frac{1}{\phi''(\psi(\beta q))}=\psi'(\beta q).$ Using definition of $\psi(\cdot),$ we obtain $\psi(\phi'(x))=x.$ Differentiating with respect to $x$ we get $\psi'(\phi'(x))\phi''(x)=1.$ This give us $\phi''(x)=\frac{1}{\psi'(\phi'(x))}.$  Let $\phi'(x)=y.$ Then $x=\psi(y).$ Evaluating  $\phi''(x)$ we get $\phi''(x)=\phi''(\psi(y)).$ Hence $\frac{1}{\psi'(y)}=\phi''(\psi(y)).$ Claim 2 proved.

\indent  From equation \eqref{eq5} it is easy to see that
\begin{equation}
g_{1}'(\beta)=\sum_{j=1}^{n-1}l_j q_j \psi'(\beta q_j)-l_n q_n \psi'(\beta q_n)
\label{eq7}
\end{equation}
where, $\frac{d\psi}{dx}\mid_{x=\beta q}:=\psi'(x).$
\indent Using Claim 2, equation \eqref{eq7} can be rewritten as \;$g_{1}'(\beta)=\sum_{j=1}^{n-1}\frac{l_j q_j}{\phi''(\psi(\beta q_j))}-\frac{l_n q_n}{\phi''(\psi(\beta q_n))}$. Hence $g'(0)<0$ implies that $g_{1}'(\beta)>0.$

\indent \textit{Sufficiency}: It easily follows from claim 2, and equation \eqref{eq7}.
\end{proof}
\end{lem}
 \indent As a consequence, we have the following result.

\begin{cor}
\label{cor2}
A critical point $x_0 $ is a  point of local maximum in side \octant \; if $g_{1}'(\beta)>0,g_{1}(\beta)=b, $  and only if $g_{1}'(\beta)\geq 0,g_{1}(\beta)=b$ for some $\beta \in \mathbb{R}^+.$

\begin{proof}
The proof easily follows from theorem \ref{thm6} and lemma \ref{lem4}.
\end{proof}
\end{cor}

 \indent Next we will show that the function $f(x)=\sum_{i=1}^{n}c_i \phi(x_i)$ has at most one point of local maximum in side \octant. Before going over the proof, we will prove some useful properties of function $g_{1}(\cdot)$, which will be used frequently in following sections.
 Notice that  $$g_{1}(\beta)=\sum_{j=1}^{n-1}l_j \psi(\beta q_j)-l_n \psi(\beta q_n)$$
 and,
 $$g_{1}'(\beta)=\sum_{j=1}^{n-1}l_j q_j \psi'(\beta q_j)-l_n q_n \psi'(\beta q_n)$$
 
 \begin{rem}
 \label{remnew}
If $q_j=q_k$ for some $j,k$ then the optimization problem can be reduced to similar problem of lower dimension. Therefore, in the future sections, we will assume that $q_j$  where $j \in \{1 \ldots n\}$ takes distinct values. 
 \end{rem}

\begin{lem}
\label{lem5}
Suppose that $q_n<q_{j0}$, where $j0 \neq n.$ Then the following statements are true. 
\begin{enumerate}[\rm(i)]
\item{$g_{1}'(\beta)$ has at most two roots on the interval $(0,\beta_{max}]$.}
\item{ $g_{1}'(\beta) \rightarrow -\infty ,$ as $\beta \rightarrow \beta_{max}.$}
\end{enumerate}

\indent \begin{proof} Denote $h(\beta, q_j,q_n):=\frac{\psi'(\beta q_j)}{\psi'(\beta q_n)},$ and $h(\beta, q_l,q_n):= \frac{\psi'(\beta q_l)}{\psi'(\beta q_n)},$ where $j,l \in \{1 \ldots n-1\}.$ Recall that $\psi'(\beta q)=\frac{d\psi}{ds}\mid_{s=\beta q}.$    

\indent \rm(i) Suppose $q_n<q_{j0} $ for some $j0 \in \{1 \ldots n-1\}.$ Let $q_1<q_2 <\ldots<q_k<q_n<q_{k+1}<\ldots<q_{n-1}.$\\  
\indent Using assumption 2, and definition of $h(\cdot,\cdot,\cdot),$ we get $h_{\beta}(\beta,q_j,q_n)\neq 0 .$ In assumption 3, we saw that  $\frac{\partial}{\partial \beta}\Big[ \frac{h_{\beta}(\beta, q_j,q_n)}{h_{\beta}(\beta, q_l,q_n)}\Big] $ has the same sign for all $(\beta,q_j,q_n,q_l)$ such that $\beta \in (0,\beta_{max}],$ and $0<q_j<q_n <q_l.$ \\ 

\indent Using assumptions 2 and 3, we obtain
\begin{align}
\frac{\partial}{\partial\beta}\Big(\log \mid \frac{h_{\beta}(\beta, q_j,q_n)}{h_{\beta}(\beta, q_l,q_n)}\mid \Big)
 &= \frac{h_{\beta}(\beta, q_l,q_n)}{h_{\beta}(\beta,q_j,q_n)}\cdot \frac{\partial}{\partial \beta} \Big(\frac{h_{\beta}(\beta,q_j,q_n)}{h_{\beta}(\beta,q_l,q_n)}\Big) \nonumber
 \neq 0. \nonumber
\end{align}

\indent Notice that for all pairs $(q_j,q_l)$ such that $q_j<q_l$ the expression $h_{\beta \beta}(\beta, q_j,q_n)h_{\beta}(\beta, q_l,q_n)-h_{\beta \beta}(\beta,q_l,q_n)h_{\beta}(\beta, q_j,q_n)$ is sign definite, where $h_{\beta \beta}(\beta, q_j,q_n)=\frac{\partial}{\partial \beta}(h_{\beta}(\beta,q_j,q_n)).$  This implies that 

$$\sum_{j=1}^{k}\sum_{l=k+1}^{n-1}\underbrace{l _j q_j l_l q_l}_{\text{is positive}}\Big[ h_{\beta \beta}(\beta, q_j,q_n)h_{\beta}(\beta, q_l,q_n)-h_{\beta}(\beta ,q_j,q_n)h_{\beta \beta}(\beta, q_l,q_n)\Big]\neq 0,$$
which, in turn implies that 
$$\sum_{j=1}^{k}\alpha_j h_{\beta \beta}(\beta, q_j,q_n)\sum_{l=k+1}^{n-1}\alpha_l h_{\beta}(\beta, q_l,q_n)-\sum_{j=1}^{k}\alpha_j h_{\beta}(\beta, q_j,q_n)\sum_{l=k+1}^{n-1}\alpha_l  h_{\beta \beta}(\beta, q_l,q_n)\neq 0$$

\indent where, $\alpha_j:=l_jq_j, $ and $\alpha_l :=l_l q_l.$ 

\indent Denote $g(\beta):=\displaystyle \frac{\sum_{j=1}^{k}\alpha_j h_{\beta}(\beta, q_j,q_n)}{\sum_{l=k+1}^{n-1}\alpha_l  h_{\beta}(\beta, q_l,q_n)}\displaystyle.$ We can see that left side of above equation  is same as numerator of $g'(\beta).$  This implies that $g'(\beta)$ is sign definite. Therefore there exists at most one $\beta$, such that $g(\beta_1)=-1.$ This implies that $\sum_{j=1}^{k}\alpha_j  h_{\beta}(\beta, q_j,q_n) +\sum_{l=k+1}^{n-1} \alpha_l  h_{\beta}(\beta, q_l,q_n)=0$ for at most one value of $\beta.$\\
\indent Using definition of $h_{\beta}(\beta, q,q_n),$ we obtain  $\sum_{j=1}^{k}\alpha_j \frac{\partial}{\partial\beta}(h(\beta, q_j,q_n))+\sum_{l=k+1}^{n-1}\alpha_l \frac{\partial}{\partial \beta}(h(\beta,q_l,q_n))-l_nq_n$  has at most one root. Recalling definition of $h(\beta,q_j,q_n)$ we obtain that function  $\sum_{j=1}^{k}\alpha_j \frac{\partial}{\partial\beta}\Big(\frac{\psi'(\beta q_j)}{\psi'(\beta q_n)}\Big)+\sum_{l=k+1}^{n-1}\alpha_l \frac{\partial}{\partial \beta}\Big(\frac{\psi'(\beta q_l)}{\psi'(\beta q_n)}\Big)-l_n q_n$  has at most one root. It is easy to see that above function is equal to  $\frac{d}{d\beta}\Big(\frac{g_{1}'(\beta)}{\psi'(\beta q_n)}\Big)$. Since $\psi'(\beta q)<0$, we obtain that $g_{1}'(\beta)$ has at most two roots in $(0,\beta_{max}].$ This completes the proof for part\rm(i).

\indent \rm(ii) Using assumption 1, $\phi '(x)$ is a decreasing function for all $x>0,$ and $\lim_{x \rightarrow \infty}\phi'(x)=0$\;(since $\lim_{y \rightarrow \infty}\phi(y)<\infty$).  Since $\psi(\cdot):=\phi'(\cdot)^{-1},$ we get $\psi:(0,\phi'(0)]\rightarrow   [0,\infty)$. Moreover, $\psi'(x)<0 $ for all $x>0.$
\indent   Using the fact that $\psi'(y)=\frac{1}{\phi''((\phi')^{-1}(y))}$( by Claim 2 in lemma \ref{lem4}) and $\phi''(0)=0,$  we obtain 
\begin{equation}
\label{eq9}
 \psi'(y) \rightarrow -\infty, \text{as} \; y \rightarrow \phi'(0).
 \end{equation}

\indent Recall that $g_{1}'(\beta)=\sum_{j=1}^{n-1}l_j q_j \psi'(\beta q_j)-l_n q_n \psi'(\beta q_n)$. Since, $q_n <q_{j0}=max(q_j)_{j=1}^{n},$ we get $q_{j0}\beta_{max}=\phi'(0).$ Hence using equation \eqref{eq9}, we obtain that $g_{1}'(\beta) \rightarrow -\infty ,$ as $\beta \rightarrow \beta_{max}.$ Proof for part \rm(ii) completed.

\end{proof}
\end{lem}

\begin{rem}
\label{rem1} 
Similarly we can show that if $q_n=\max(q_j)_{j=1}^{n},$ then $g_{1}'(\beta) \rightarrow \infty ,$ as $\beta \rightarrow \beta_{max}.$
\end{rem}

 \indent \begin{lem}
 \label{lem6}

 Suppose that $q_n >q_j$ for all $j \in \{1\ldots n-1\}.$ Then, $g_{1}'(\beta)$ has at most one root.  Moreover, if $\lim_{\beta \rightarrow 0}g_{1}'(\beta)\geq 0$, then $g_{1}'(\beta)>0$ for all $\beta \in (0, \beta_{max}].$\\
 \indent \begin{proof} Assume $q_n>q_j$ for all $j \in \{1\ldots n-1\}.$ We will use contrapositive approach to prove this claim. Suppose there exist $\beta_1 \neq \beta_2$ such that $g_{1}'(\beta_1)=g_{1}'(\beta_2)=0.$ Since  $g_{1}'(\beta_1)=0$, we get $\displaystyle \frac{\sum_{j=1}^{n-1}l_j q_j \psi'(\beta_1 q_j)}{l_n q_n \psi'(\beta_1 q_n)}=1.$  Similarly, $\displaystyle \frac{\sum_{j=1}^{n-1}l_j q_j \psi'(\beta_2 q_j)}{l_n q_n \psi'(\beta_2 q_n)}=1.$
\indent Denote $\frac{l_j q_j}{l_n q_n}:=\alpha_j > 0 $ for all $j$, and   $F(\beta):=\sum_{j=1}^{n-1}\alpha_j \frac{ \psi'(\beta q_j)}{ \psi'(\beta q_n)}.$ Hence  there exist $\beta_1 \neq \beta_2$ such that $F(\beta_1)=F(\beta_2)=1.$ So there exist $\beta_0$ such that $F'(\beta_0)=0.$ This implies that 
$\sum_{j=1}^{n-1}\alpha_j \frac{d}{d\beta}\Big(\frac{\psi'(\beta q_j)}{\psi'(\beta q_n)}\Big)\mid_{\beta=\beta_0}=0,$
where, $q_j <q_n$ for all $j \in \{1 \ldots n-1\}.$
\indent Since $\alpha_j > 0$ for all $j \in \{1 \ldots n-1\}$, we conclude that $\frac{d}{d\beta}\Big(\frac{\psi'(\beta q_j)}{\psi'(\beta q_n)}\Big)\vert_{\beta=\beta_0}$ is not sign definite. Hence we obtain contradiction to assumption 2,  which says that   $\frac{d}{d\beta}\Big(\frac{\psi'(\beta p)}{\psi'(\beta q)}\Big),$ where $p \neq q$ is sign definite. Therefore, if $q_n > q_j$ for all $j \in \{1 \ldots n-1\},$ then $g_{1}'(\beta)$ has at most one root in the interval $(0, \beta_{max})$.

\indent Moreover, assume that $\lim_{\beta \rightarrow 0}g_{1}'(\beta)\geq 0$, where  $g_{1}'(\beta)=\sum_{j=1}^{n-1}l_j q_j \psi'(\beta q_j)-l_n q_n \psi'(\beta q_n).$ First, consider the case when $\lim_{\beta \rightarrow 0}g_{1}'(\beta)> 0$. Then, the function $g_{1}(\cdot)$ can either be identically increasing or $g_{1}(\cdot)$ changes monotonicty. We will show that the second case cannot happen. Using the remark \ref{rem1} we obtain
 
 \begin{equation}
 \label{eq10}
 g_{1}'(\beta)\rightarrow \infty \text{ as}\; \beta \rightarrow \beta_{max}.
 \end{equation}

\indent Suppose by contradiction that $g_{1}(\cdot)$ changes monotonicity. Using the fact $\lim_{\beta \rightarrow 0}g_{1}'(\beta)>0$, and \eqref{eq10}, we get  $g_{1}(\beta)$ has two critical points, given by $\overline{\beta}$ and $\beta_1.$ But this contradicts first part of this lemma. Hence $g_{1}'(\beta)>0$ for all $\beta.  \in (0,\beta_{max}].$\\
Next let $\lim_{\beta \rightarrow 0}g_{1}'(\beta)= 0$. Using remark \ref{rem1} and first part of this lemma, we obtain that  $g_{1}'(\beta)>0 $ for all $\beta \in (0,\beta_{max}].$
\end{proof}
\end{lem}

\begin{rem}
\label{rem2}
In lemmas \ref{lem5} and \ref{lem6}, we assumed that the last component is negative. Similar results hold true  for side orthant with first or second component negative. These results will be used in the proof for case of two side orthants. 
\end{rem}

 \indent Next, we will show the main result of the section.

\begin{prop}
\label{prop8}
The function $f(x)= \sum_{i=1}^{n} c_i \phi(x_i)$ has at most one  point of local maximum in the side \octant.
 \begin{proof} Recall that, $g_{1}(\beta)=\sum_{j=1}^{n-1}l_j \psi(\beta q_j)-l_n \psi(\beta q_n),$ and $g_{1}'(\beta)=\sum_{j=1}^{n-1}l_j q_j \psi'(\beta q_j)-l_n q_n \psi'(\beta q_n).$
  \indent \underline{Case 1}:  Suppose $q_{n}> q_{j}$ for all $j \in \{ 1 \ldots n-1\}.$  We need to show that for a given $b \in \mathbb{R}$, there exists at most one value of  $\beta \in (0,\beta_{max})$ such that $g_{1}(\beta)=b, $ and $g_{1}'(\beta)>0$(corollary \ref{cor2}). Using remark \ref{rem1}, we get that $g_{1}'(\beta) \rightarrow \infty, $ as $\beta \rightarrow \beta_{max}.$ Moreover, from  lemma \ref{lem6},we obtain that  that if $q_n>q_j $ for all $j \in \{ 1\ldots n-1\},$ then $g_{1}'(\beta)$ has at most one root. Hence, there exists at most one value of $\beta$ such that $g_{1}(\beta)=b, $ and $g_{1}'(\beta)>0.$ This completes the proof for Case 1.

\indent \underline{Case 2}: Suppose $q_n<q_{j0}$ where $j0 \neq n.$
 Using lemma \ref{lem5} we get that $g_{1}'(\beta)$ has at most two roots in the interval $(0,\beta_{max}].$ In addition, we showed in lemma \ref{lem5}  that $ g_{1}'(\beta) \rightarrow -\infty ,$ as $\beta \rightarrow \beta_{max}.$ Hence, for a given $b \in \mathbb{R},$ there exists at most one $\beta \in (0,\beta_{max}]$ such that $g_{1}(\beta)=b,$ and $g_{1}'(\beta)\geq 0.$ Therefore $f(x)$ has at most one point of local maximum in side \octant. Case 2 completed.
\end{proof}
\end{prop}

\indent In  section 4, we saw that the function $f(x)$ has at most one point of local maximum in main \octant \; or side \octant\;  with one component negative. It might happen that there are two points of local maxima,  one in main \octant \; and other in side \octant \; with one component negative.  But in the next section we will show that $f(x)$ does not have points of local maxima in both main \octant \; and side \octant \; with one negative component.

\section{ Points of Local Maxima in Main and Side Orthant}

\indent In this section we will show that the function $f(x)=\sum_{i=1}^{n} c_i \phi(x_i)$, defined on the hyperplane P, does not have points of local maxima in both main \octant \; and side \octant. In the main \octant \;  $f(x)$ is of the form 
$$f_{1}(\beta)=\sum_{i=1}^{n}l_i \psi(\beta q_j),$$ and in side \octant \; it is of the form 
$$g_{1}(\beta)=\sum_{i=1}^{n-1}l_i \psi(\beta q_i)-l_n \psi(\beta q_n).$$

 \indent We will first present two auxiliary  results. These will be used to show that $f(x)$ does not have points of local maxima in both main \octant \; and side \octant.
 \begin{lem}
 \label{lem7}
 Suppose $f_{1}(\beta)$ and $g_{1}(\beta)$ are defined as above.  Then the following are true:
 \begin{enumerate}[\rm(i)]
 \item{$ f_{1}(\beta) \geq g_{1}(\beta)$ for all $\beta \in (0, \beta_{max}]$}
 \item{$f_{1}'(\beta)<0$ for all $\beta 
 \in (0, \beta_{max}].$}
 \end{enumerate}
 \indent \begin{proof}
 It can be easily checked that 
 \begin{equation}
 \label{eq11}f_{1}'(\beta)=\sum_{j=1}^{n} l_jq_j\psi'(\beta q_j),
 \end{equation}
 and, 
 \begin{equation}
 \label{eq12}
  g_{1}'(\beta)=\sum_{j=1}^{n-1}l_j q_j \psi'(\beta q_j)-l_n q_n \psi'(\beta q_n).
  \end{equation}
  where $\psi(\cdot)=(\phi'(\cdot))^{-1},\psi'(\beta q)=\frac{d \psi}{ds}\mid_{s=\beta q}.$\\
  \begin{enumerate}[(i)]
 \item{  We will first show that $f_{1}(\beta) \geq g_{1}(\beta)$ for all $\beta \in (0, \beta_{max}]$. 
 \indent Using assumption  1, we have $\phi'(x)\geq 0$ for all $x.$ Hence $\psi(x)\geq 0 $ for all $x.$ Moreover,  it is easy to see that the expressions for $f_{1}(\beta)$ and $g_{1}(\beta)$ are identical except the term $l_n  \psi(\beta q_n).$ Since $l_n \psi(\beta q_n)\geq0$, we obtain $  f_{1}(\beta) \geq g_{1}(\beta)$ for all $\beta \in (0, \beta_{max}].$ This completes the proof for (i).}
 
 \item{ The proof follows from the fact that $\psi'(x)<0$ for all $x$\; ( assumption 1 ).}
 \end{enumerate}
 \end{proof}
 \end{lem}

\indent  Now, we will present the main result of this section.

\begin{prop}
\label{prop9}
The function $f(x)=\sum_{i=1}^{n} c_i \phi(x_i)$ has at most one point of local maximum in  main \octant \; and side \octant .\\
\indent \begin{proof} We need to show that $f(x)$ does not have points of local maxima in both main \octant \; and side \octant. Using proposition \ref{cor1} and corollary \ref{cor2}, we need to show that for given $b \in \mathbb{R},$ there does not exist $\beta_1 \neq \beta_2 \in (0,\beta_{max}]$ such that $f_{1}(\beta_1)=b,$ and $g_{1}(\beta_2)=b, g_{1}'(\beta_2)\geq0.$  We will prove the above result using two cases.

\indent \underline{Case 1}: Let $q_n>q_j$ for all $j \in \{1 \ldots n-1\}.$ Using assumption 1, we get that $\psi'(y)<0,$ and $\psi(y) \in [0,\infty)$ where $y\in (0,\phi'(0)].$ We obtain that $\lim_{y \rightarrow \phi'(0)}\psi(y)=0.$  Define $y:=\beta q, $ then $\lim_{\beta \rightarrow \beta_{max}}\psi(\beta q_n)=0, $ which in turn gives $\lim_{\beta \rightarrow \beta _{max}}l_n \psi(\beta q_n)=0.$ Hence 
\begin{equation}
\label{eq13}
f_{1}(\beta_{max)})=g_{1}(\beta_{max}).
\end{equation}

\noindent Now we will look at following possibilities.

\indent \underline{Sub Case \rm(i)}: Suppose $\lim_{\beta \rightarrow 0}g_{1}'(\beta)<0.$ First, we will consider the case when $g_{1}(\beta)$ is decreasing identically. Using corollary \ref{cor2}, $f(x)$  does not have any point of local maximum in side \octant. So $f(x)$ has at most one point of local maximum in the main \octant.\\
\indent Next, suppose that $g_{1}(\beta)$ changes monotonicity. We showed in lemma\; \ref{lem6} that $g_{1}'(\beta)$ has at most one root.  Denote $\overline{\beta} \in [0,\beta_{max}]$ as the point where $g_{1}(\beta)$ changes monotonicity. To summarize, $f_{1}(\beta)$ and $g_{1}(\beta)$ satisfy the following conditions :
\begin{enumerate}[\rm(i)]
\item{$f_{1}'(\beta)<0,$ and $ f_{1}(\beta)> g_{1}(\beta)$ for all $\beta \in (0,\beta_{max})$; see lemma \ref{lem7}}
\item{$g_{1}'(\beta)<0$ on the interval $(0,\overline{\beta}]$ and $g_{1}'(\beta)>0$ on the interval $[\overline{\beta},\beta_{max}].$}
\item{ $f_{1}(\beta_{max})=g_{1}(\beta_{max}).$}
\end{enumerate}

\indent Using the above conditions, we can deduce that for a given $b \in \mathbb{R}$ there does not exist $\beta_1 \neq \beta_2$ such that $f_{1}(\beta)=b $  and $g_{1}(\beta)=b,g_{1}'(\beta)>0.$  Hence, $f(\cdot)$ has at most one point of local maximum in either main or side \octant. This completes the proof for Sub Case \rm(i).

\indent \underline{Sub Case \rm(ii)}: Suppose $\lim_{\beta \rightarrow 0}g_{1}'(\beta)\geq 0.$ From lemma \ref{lem6}, we get that  $g_{1}(\beta)$ is a monotonically increasing function. Using lemma \ref{lem7}, $f_{1}(\beta)> g_{1}(\beta) $ for all $\beta \in [0,\beta_{max})$, and $f_{1}'(\beta)<0$ for all $\beta.$  At the same time, $g_{1}(\beta_{max})=f_{1}(\beta_{max}).$  Hence, for a given $b \in \mathbb{R}$,  there exists at most one $\beta \in (0,\beta_{max})$ such that either $f_{1}(\beta)=b, $ or $g_{1}(\beta)=b, g_{1}'(\beta)>0.$\\ Proof for Sub Case 2 completed.
\indent Combining the results of Sub Case 1 and Sub Case 2, we conclude that  $f(x)$ does not have points of local maxima in both main \octant \; and side \octant. This completes the proof for Case 1.

\indent \underline{ Case 2}: Suppose $q_n<q_{j0},$ where $j0 \in \{1 \ldots n-1\}.$  We will go over this case by contradiction. Suppose the function $f(x)$ has two points of local maxima, one in each main and side \octant. This implies that there exist $\beta_1 \neq \beta_2, $ such that for a given $b \in \mathbb{R},$ \; $f_{1}(\beta_1)=g_{1}(\beta_2)=b, $ and $g_{1}'(\beta_2)\geq 0.$ Notice that $\beta_1 $ and $\beta_2$ lie in $ (0,\beta_{max}].$  Recall that $f_{1}(\beta)=\sum_{j=1}^{n}l_j \psi(\beta q_j)$, and
$g_{1}(\beta)=\sum_{j=1}^{n-1}l_j \psi(\beta q_j)-l_n \psi(\beta q_n).$

\noindent Using part \rm (ii) of lemma \ref{lem5}, we get $g_{1}'(\beta) \rightarrow -\infty,$ as $\beta \rightarrow \beta_{max}.$ In addition, we assume that $g_{1}'(\beta_2) >0.$ This implies that  there exists $\beta' \in (\beta_2, \beta_{max})$ such that $g_{1}'(\beta')=0,$ and $g_{1}(\beta')>b=f_{1}(\beta_1)> f_{1}(\beta_{max})$( lemma \ref{lem7}).\\
\indent Since $g_{1}'(\beta')=0,$ we get $l_n=\frac{1}{q_n \psi'(\beta' q_n)}\sum_{j=1}^{n-1}l_j q_j \psi'(\beta' q_j).$ This implies that 
$g_{1}(\beta')= \sum_{j=1}^{n-1}l_j\Big(\psi(\beta' q_j)-\frac{\psi(\beta' q_n)}{q_n \psi'(\beta' q_n)}q_j \psi'(\beta' q_j)\Big)$. Moreover, we also have 
$f_{1}(\beta)=\sum_{j=1}^{n-1}l_j \Big(\psi(\beta q_j)+\frac{\psi(\beta q_n)}{q_n \psi'(\beta' q_n)}q_j \psi'(\beta' q_j)\Big)$. Then, $f_{1}(\beta_{max})=\sum_{j=1}^{n-1}l_j \Big(\psi(\beta_{max} q_j)+\frac{\psi(\beta_{max} q_n)}{q_n \psi'(\beta_{max} q_n)}q_j \psi'(\beta_{max} q_j)\Big).$ Since $g_{1}(\beta') > f_{1}(\beta_{max}),$ we obtain 
$$\sum_{j=1}^{n-1}l_j\Big(\psi(\beta' q_j)-\frac{\psi(\beta' q_n)}{q_n \psi'(\beta' q_n)}q_j \psi'(\beta' q_j)\Big)>\sum_{j=1}^{n-1}l_j \Big(\psi(\beta_{max} q_j)+\frac{\psi(\beta_{max} q_n)}{q_n \psi'(\beta_{max} q_n)}q_j \psi'(\beta_{max} q_j)\Big).$$ We saw in section 4, that $l_j>0$ for all $j \in \{1 \ldots n\}.$ Hence there exists $j \in \{1 \ldots n-1\}$ such that 
$$\psi(\beta' q_j)-\frac{\psi(\beta' q_n)}{q_n \psi'(\beta' q_n)}q_j \psi'(\beta' q_j)>\psi(\beta_{max} q_j)+\frac{\psi(\beta_{max} q_n)}{q_n \psi'(\beta_{max} q_n)}q_j \psi'(\beta_{max} q_j).$$

\indent We will show that the above inequality is not true. In other words, we will show that for all $q_j,$ where $j \in \{1 \ldots n-1\},$
\begin{equation}
\label{eq14}
\psi(\beta' q_j)-\frac{\psi(\beta' q_n)}{q_n \psi'(\beta' q_n)}q_j \psi'(\beta' q_j)\leq\psi(\beta_{max} q_j)+\frac{\psi(\beta_{max} q_n)}{q_n \psi'(\beta_{max} q_n)}q_j \psi'(\beta_{max} q_j).
\end{equation}

\noindent  Pick arbitrary $q_j,$ where $j \in \{1 \ldots n-1\}.$ Then we will consider the following possibilities. \\

\indent \underline{ Sub Case \rm(i)}: Let $q_n >q_j$. Using assumption 2, we have $\frac{\psi'(\beta p)}{\psi '(\beta q)}, p>q$ is an increasing function of $\beta.$ This implies $-\frac{q_j}{q_n} \psi(\beta q_n)\frac{d}{d\beta}\Big(\frac{\psi'(\beta q_j)}{\psi '(\beta q_n)}\Big)>0$(since $\psi(\cdot)>0$). Notice that 
\begin{align}
\frac{d}{d\beta}\Big(\psi(\beta' q_j)-\frac{\psi(\beta' q_n)}{q_n \psi'(\beta' q_n)}q_j \psi'(\beta' q_j)\Big)&=q_j \psi'(\beta' q_j)-\frac{q_j}{q_n} \psi(\beta' q_n) \frac{d}{d\beta}\Big(\frac{\psi'(\beta' q_j)}{\psi'(\beta' q_n)}\Big)-\frac{q_j}{q_n}\frac{\psi'(\beta' q_j)}{\psi'(\beta' q_n)}q_n \psi'(\beta' q_n)\nonumber\\
&=-\frac{q_j}{q_n} \psi(\beta' q_n)\frac{d}{d\beta}\Big(\frac{\psi'(\beta' q_j)}{\psi '(\beta' q_n)}\Big)\nonumber\\
&>0.\nonumber
\end{align}

\noindent The left hand side of inequality \eqref{eq14} is an increasing function of $\beta$. Hence, we obtain  
\begin{align}
\psi(\beta' q_j)-\frac{\psi(\beta' q_n)}{q_n \psi'(\beta' q_n)}q_j \psi'(\beta' q_j)&<\psi(\beta_{max} q_j)-\frac{\psi(\beta_{max} q_n)}{q_n \psi'(\beta_{max} q_n)}q_j \psi'(\beta_{max} q_j)\nonumber\\
&<\psi(\beta_{max} q_j)+\frac{\psi(\beta_{max} q_n)}{q_n \psi'(\beta_{max} q_n)}q_j \psi'(\beta_{max} q_j).
\end{align}

since $\displaystyle \frac{\psi(\beta_{max} q_n)}{q_n \psi'(\beta_{max} q_n)}q_j \psi'(\beta_{max} q_j)>0.$\\
\noindent Therefore inequality \eqref{eq14} is true. This completes the proof for Sub Case \rm(i).\\
\indent \underline{Sub Case \rm(ii)}: Suppose $q_n <q_j.$ Using assumption 4 we have $\frac{d}{d\beta}\Big(\frac{\psi(\beta q)}{\psi(\beta p)} \Big)>0,$ where $p>q.$ Then we obtain a sequence of inequalities,

\begin{align}
& \frac{d}{d\beta}\Big(\ln  \frac{\psi(\beta q)}{\psi(\beta p)} \Big)>0,\nonumber\\
& \frac{d}{d\beta}\Big(\ln (\psi(\beta q))\Big)>\frac{d}{d\beta}\Big(\ln (\psi(\beta p))\Big)\nonumber\\
&\frac{1}{\psi(\beta p)}\frac{d}{d\beta}\Big(\psi(\beta p)\Big)<\frac{1}{\psi(\beta q)}\frac{d}{d\beta}\Big(\psi(\beta q)\Big)\nonumber\\
&\psi(\beta p)<\psi(\beta  q)\frac{\frac{d}{d\beta}(\psi(\beta p))}{\frac{d}{d\beta}(\psi(\beta q))}.\nonumber
\end{align}
since  $\psi'(x)<0,$ and $ \psi(x)>0.$\\

\noindent Let $p:=q_j,q:=q_n, \beta:=\beta'$, we obtain $\psi(\beta' q_j)-\frac{\psi(\beta' q_n)}{q_n \psi'(\beta' q_n)}q_j \psi'(\beta' q_j)<0.$ Hence left hand side of inequality \eqref{eq14} is negative.\\
\noindent We have seen earlier that $\psi(x)\geq 0 $ for all $x \in (0,\phi'(0)],$ and $\psi'(x)<0 $ for all $x \in (0, \phi'(0)).$ This implies that $\psi(\beta_{max} q_j)+\frac{\psi(\beta_{max} q_n)}{q_n \psi'(\beta_{max} q_n)}q_j \psi'(\beta_{max} q_j)>0.$ Therefore inequality  \eqref{eq14} is true. This completes the proof for Sub  Case \rm(ii).\\
\indent \underline{Sub Case \rm(iii)}: Assume that $q_n=q_j.$ Then it is easy to see that $\psi(\beta' q_j)-\frac{\psi(\beta' q_n)}{q_n \psi'(\beta' q_n)}q_j \psi'(\beta' q_j)=0.$ The left hand side of inequality \eqref{eq14} is zero.
The right hand side of inequality \eqref{eq14} is equal to $2\psi(\beta_{max}q_n).$ If $\beta_{max} q_n=\phi'(0),$ then $\psi(\beta_{max} q_n)=0.$ Hence we get equality.\\
\indent Therefore, we obtain contradiction to our assumption. Combining the results of Case 1 and Case 2, we obtain that $f(x)$ does not have points of local maxima in main \octant \; and side \octant. 
\end{proof}
\end{prop}

\indent We have seen that the function $f(x)=\sum_{i=1}^{n} c_i \phi(x_i), c_i \neq 0$ for all $i \in \{1 \ldots n\}$ does not have points of local maxima in main and side \octant. Now, we will show similar result for two side \octants.\\

\section{ Points of local maxima in Two Side Orthants}

In this section, we will show that function $f(x)=\sum_{i=1}^{n}c_i \phi(x_i), c_i \neq 0$ for all $n$, does not have points of local maxima in two side \octants. We will show the proof for the case, when the first and second components are negative. Other cases can be analyzed similarly. \\

\noindent In the first side \octant,\; $f(x)$ takes the form 
\begin{equation}
\label{eq15}
g_{1}(\beta)=\sum_{j=2}^{n}l_j \psi(\beta q_j)-l_1 \psi(\beta q_1),
\end{equation}
and in the second side \octant \; $f(x)$ is of the form
\begin{equation}
\label{eq16}
g_{2}(\beta)=\sum_{j=1,j \neq 2}^{n}l_j \psi(\beta q_j)-l_{2}\psi(\beta q_2).
\end{equation}

\noindent It can be easily checked that 
\begin{equation}
\label{eq17}
g_{1}'(\beta)=\sum_{j=2}^{n}l_j q_j \psi'(\beta q_j)-l_1 q_1 \psi'(\beta q_1),
\end{equation}
and 
\begin{equation}
\label{eq18}
g_{2}'(\beta)=\sum_{j=1,j \neq 2}^{n}l_j q_j \psi'(\beta q_j)-l_2 q_2 \psi'(\beta q_2).
\end{equation}

\indent First, we will solve an auxiliary problem. The main result of this section will be an easy consequence of the solution to this auxiliary problem. 

\subsection{Three Point Problem}
\begin{lem}
\indent \label{lem8}
Consider the functions $g_{1}(\cdot)$ and $g_{2}(\cdot)$ defined above. Suppose $\beta_1 ,$ and $\beta_2$ are critical points for $g_{1}(\cdot)$ and $g_{2}(\cdot)$ respectively. Moreover assume that $g_{1}(\beta_1)-g_{2}(\beta_2)<0.$ Then the following inequality holds:

 \begin{align}
 &\frac{\beta_2}{2\beta_{1}^{2}}\Big[\frac{\psi_{\beta}(\beta_2 q_3)}{\psi_{\beta}(\beta_1 q_3)}\Big(\frac{\psi(\beta_1 q_1)}{\psi_{\beta}(\beta_1 q_1)}-\frac{\psi(\beta_1 q_2)}{\psi_{\beta}(\beta_1 q_2)}\Big)\\
 &+\frac{\psi_{\beta}(\beta_2 q_3)}{\psi_{\beta}(\beta_1 q_3)}\cdot \frac{\psi_{\beta}(\beta_2 q_1)}{\psi_{\beta}(\beta_1 q_1)}\Big(\frac{\psi(\beta_2 q_1)}{\psi_{\beta}(\beta_2 q_1)}-\frac{\psi(\beta_2 q_2)}{\psi_{\beta}(\beta_2 q_2)}\Big)\\
 &+\frac{\psi_{\beta}(\beta_2 q_2)}{\psi_{\beta}(\beta_1 q_2)}\Big(\frac{\psi(\beta_1 q_1)}{\psi_{\beta}(\beta_1 q_1)}-\frac{\psi(\beta_1 q_3)}{\psi_{\beta}(\beta_1 q_3)}\Big)\\
 &+\frac{\psi_{\beta}(\beta_2 q_2)}{\psi_{\beta}(\beta_1 q_2)}\cdot \frac{\psi_{\beta}(\beta_2 q_1)}{\psi_{\beta}(\beta_1 q_1)}\Big(\frac{\psi(\beta_2 q_1)}{\psi_{\beta}(\beta_2 q_1)}-\frac{\psi(\beta_2 q_3)}{\psi_{\beta}(\beta_2 q_3)}\Big)\\
 &+\frac{\psi_{\beta}(\beta_2 q_1)}{\psi_{\beta}(\beta_1 q_1)}\Big(\frac{\psi(\beta_1 q_3)}{\psi_{\beta}(\beta_1 q_3)}-\frac{\psi(\beta_1 q_2)}{\psi_{\beta}(\beta_1 q_2)}\Big)\\
 &+\Big(\frac{\psi(\beta_2 q_2)}{\psi_{\beta}(\beta_2 q_2)}-\frac{\psi(\beta_2 q_3)}{\psi_{\beta}(\beta_2 q_3)}\Big)\Big(\frac{\psi_{\beta}(\beta_2 q_3)}{\psi_{\beta}(\beta_1 q_3)}\cdot \frac{\psi_{\beta}(\beta_2 q_1)}{\psi_{\beta}(\beta_1 q_1)}-\frac{\psi_{\beta}(\beta_2 q_1)}{\psi_{\beta}(\beta_1 q_1)}\cdot\frac{\psi_{\beta}(\beta_2 q_2)}{\psi_{\beta}(\beta_1 q_2)}-\frac{\psi_{\beta}(\beta_2 q_3)}{\psi_{\beta}(\beta_1 q_3)}\cdot \frac{\psi_{\beta}(\beta_2 q_2)}{\psi_{\beta}(\beta_1 q_2)}\Big)\Big]\nonumber\\
 &<0.\nonumber
 \end{align}
Here, $\psi_{\beta}(\beta_i  q_j):=\frac{d}{d\beta}(\psi(\beta q_j))\vert_{\beta=\beta_i}$, and $i \in\{1,2\},j \in \{1,2,3\}.$ \begin{proof} Notice that the proof is quite technical. Since $g_{1}'(\beta_1)=0,$ and $g_{2}'(\beta_2)=0$, we obtain
$$l_1 q_1 \psi'(\beta_1 q_1)=l_2 q_2 \psi'(\beta_1 q_2)+\sum_{j=3}^{n}l_j q_j \psi'(\beta_1 q_j),\; \text{and}$$ 
$$l_2 q_2 \psi'(\beta_2 q_2)=l_1 q_1 \psi'(\beta_2 q_1)+\sum_{j=3}^{n}l_j q_j \psi'(\beta_2 q_j)$$

\indent Solving for $l_1$ and $l_2$ we get

\begin{equation}
\label{eq19}
l_1=\frac{1}{q_1}\sum_{j=3}^{n}l_j q_j \Big(\frac{\psi'(\beta_1 q_2)\psi'(\beta_2 q_j)+\psi'(\beta_2 q_2)\psi'(\beta_1 q_j)}{\psi'(\beta_1 q_1)\psi'(\beta_2 q_2)-\psi'(\beta_2 q_1)\psi'(\beta_1 q_2)}\Big), \text{and}
\end{equation}
\begin{equation}
\label{eq20}
l_2=\frac{1}{q_2}\sum_{j=3}^{n}l_j q_j \Big(\frac{\psi'(\beta_2 q_1)\psi'(\beta_1 q_j)+\psi'(\beta_1 q_1)\psi'(\beta_2 q_j)}{\psi'(\beta_1 q_1)\psi'(\beta_2 q_2)-\psi'(\beta_2 q_1)\psi'(\beta_1 q_2)}\Big)
\end{equation}
 \noindent For simplicity, denote $\psi'(\beta_1 q_1)\psi'(\beta_2 q_2)-\psi'(\beta_2 q_1)\psi'(\beta_1 q_2):=\Delta.$ Since $l_j>0$ for all $j,$ and $\psi'(\cdot)<0,$  we get $\Delta>0.$\\

 \indent Using equations \eqref{eq19} and \eqref{eq20} we can rewrite the expressions for $g_{1}(\beta_1)$ and $g_{2}(\beta_2)$ as

 \begin{align}
  g_{1}(\beta_1)&=\sum_{j=3}^{n}l_j\Big[q_j \frac{\psi(\beta_1 q_2)}{q_2}\Big(\frac{\psi'(\beta_2 q_1)\psi'(\beta_1 q_j)+\psi'(\beta_2 q_j)\psi'(\beta_1q_1)}{\Delta}\Big)
 - \frac{q_j \psi(\beta_1 q_1)}{q_1}\Big(\frac{\psi'(\beta_1 q_2)\psi'(\beta_2 q_j)+\psi'(\beta_2 q_2)\psi'(\beta_1 q_j)}{\Delta}\Big)\nonumber\\
 &+\psi(\beta_1 q_j)\Big], \text{and}
 \end{align}
 
 \begin{align}
 g_{2}(\beta_2)&=\sum_{j=3}^{n}l_j\Big[q_j \frac{\psi(\beta_2 q_1)}{q_1}\Big(\frac{\psi'(\beta_1 q_2)\psi'(\beta_2 q_j)+\psi'(\beta_1 q_j)\psi'(\beta_2 q_2)}{\Delta}\Big)
 - \frac{q_j \psi(\beta_2 q_2)}{q_2}\Big(\frac{\psi'(\beta_2 q_1)\psi'(\beta_1 q_j)+\psi'(\beta_1 q_1)\psi'(\beta_2 q_j)}{\Delta}\Big)\nonumber\\
 &+\psi(\beta_2 q_j)\Big]
 \end{align}
 \noindent Since $g_{1}(\beta_1)-g_{2}(\beta_2)<0$, we get

 \begin{align}
 g_{1}(\beta_1)-g_{2}(\beta_2)=\nonumber\\
 &\sum_{j=3}^{n}l_j\Big[-\Big(\frac{\psi'(\beta_1 q_2)\psi'(\beta_2 q_j)+\psi'(\beta_2 q_2)\psi'(\beta_1q_j)}{\Delta}\Big)\Big(\frac{q_j \psi(\beta_1 q_1)}{q_1}+\frac{q_j \psi(\beta_2 q_1)}{q_1}\Big)\nonumber\\
 &+\Big(\frac{\psi'(\beta_2 q_1)\psi'(\beta_1 q_j)+\psi'(\beta_1 q_1)\psi'(\beta_2 q_j)}{\Delta}\Big)\Big(\frac{q_j \psi(\beta_1 q_2)}{q_2}+\frac{q_j \psi(\beta_2 q_2)}{q_2}\Big)\nonumber\\
 &+(\psi(\beta_1 q_j)-\psi(\beta_2 q_j))\Big]<0.\nonumber
 \end{align}
 \noindent Multiplying both sides of the inequality by $4\beta_1 \beta_2 q_1 q_2 \Delta$ we get 
 \begin{align}
 &\sum_{j=3}^{n}l_j\Big[-\Big(4\beta_1 \beta_2 q_2 q_j \psi'(\beta_1 q_2)\psi'(\beta_2 q_j) +4\beta_1 \beta_2 q_2 q_j \psi'(\beta_2 q_2)\psi'(\beta_1 q_j)\Big)(\psi(\beta_1 q_1)+\psi(\beta_2 q_1))\nonumber\\
 &+\Big(4\beta_1 \beta_2 q_1 q_j \psi'(\beta_2 q_1)\psi'(\beta_1 q_j)+4 \beta_1 \beta_2 q_1 q_j \psi'(\beta_2 q_j)\psi'(\beta_1 q_1)\Big)(\psi(\beta_1 q_2)+\psi(\beta_2 q_2))\nonumber\\
 &+4\beta_1 \beta_2 q_1 q_2 \Delta(\psi(\beta_1 q_j)- \psi(\beta_2 q_j))\Big]<0.\nonumber
 \end{align}
 
 \indent Since $l_j >0 $ for all $j, $ at least one of the coefficients  in the above sum should be negative. Without loss of generality let $j=3,$ and we obtain
 
 \begin{align}
&-4\beta_1 \beta_2 q_2 q_3 \Big[\psi'(\beta_1 q_2)\psi'(\beta_2 q_3)+\psi'(\beta_2 q_2)\psi'(\beta_1 q_3)\Big)\Big(\psi(\beta_1 q_1)+\psi(\beta_2 q_1)\Big)\nonumber\\
&+4\beta_1 \beta_2 q_1 q_3 \Big(\psi'(\beta_2q_1) \psi'(\beta_1 q_3)+\psi'(\beta_2 q_3)\psi'(\beta_1 q_1)\Big)\Big(\psi(\beta_1 q_2)+\psi(\beta_2 q_2)\Big)\nonumber\\
&+4\beta_1 \beta_2 q_1 q_2 \Delta \Big(\psi(\beta_1 q_3)-\psi(\beta_2 q_3)\Big)\Big]<0.\nonumber
\end{align}

\noindent We saw earlier that $\psi'(x)<0.$  Dividing both sides of above inequality by $-8\beta_{1}^{3}q_1 q_2 q_3 \psi'(\beta_1 q_1)\psi'(\beta_1 q_2)\psi'(\beta_1 q_3)$, we get 
\begin{align}
&\frac{1}{2\beta_1 q_1 \psi'(\beta_1 q_1)}\Big(\frac{2\beta_2 q_2 \psi'(\beta_2 q_2)}{2 \beta_1 q_2 \psi'(\beta_1 q_2)}+\frac{2 \beta_2 q_3 \psi'(\beta_2 q_3)}{2 \beta_1 q_3 \psi'(\beta_1 q_3)}\Big)(\psi(\beta_1 q_1)+\psi(\beta_2 q_1))\nonumber\\
&-\frac{1}{2\beta_1 q_2 \psi'(\beta_1 q_2)}\Big(\frac{2\beta_2 q_1 \psi'(\beta_2 q_1)}{2 \beta_1 q_1 \psi'(\beta_1 q_1)}+\frac{2 \beta_2 q_3 \psi'(\beta_2 q_3)}{2 \beta_1 q_3 \psi'(\beta_1 q_3)}\Big)(\psi(\beta_1 q_2)+\psi(\beta_2 q_2))\nonumber\\
&-\frac{1}{2\beta_1 q_3 \psi'(\beta_1 q_3)}\Big(\frac{2\beta_2 q_2 \psi'(\beta_2 q_2)}{2 \beta_1 q_2 \psi'(\beta_1 q_2)}-\frac{2 \beta_2 q_1 \psi'(\beta_2 q_1)}{2 \beta_1 q_1 \psi'(\beta_1 q_1)}\Big)(\psi(\beta_1 q_3)-\psi(\beta_2 q_3))<0.\nonumber
\end{align}

\indent  After distributing the terms on left side we obtain 
\begin{align}
\label{equuuu}
&\frac{\beta_2}{2\beta_{1}^{2}}\Big[\Big(\frac{\psi_{\beta}(\beta_2 q_2)}{\psi_{\beta}(\beta_1 q_2)}+\frac{\psi_{\beta}(\beta_2 q_3)}{\psi_{\beta}(\beta_1 q_3)}\Big)\Big(\frac{\psi(\beta_1 q_1)}{\psi_{\beta}(\beta_1 q_1)}+\frac{\psi(\beta_2 q_1)}{\psi_{\beta}(\beta_1 q_1)}\Big)\nonumber\\
&+\Big(\frac{\psi_{\beta}(\beta_2 q_1)}{\psi_{\beta}(\beta_1 q_1)}+\frac{\psi_{\beta}(\beta_2 q_3)}{\psi_{\beta}(\beta_1 q_3)}\Big)\Big(-\frac{\psi(\beta_1 q_2)}{\psi_{\beta}(\beta_1 q_2)}-\frac{\psi(\beta_2 q_2)}{\psi_{\beta}(\beta_1 q_2)}\Big)\nonumber\\
&+\Big(\frac{\psi_{\beta}(\beta_2 q_2)}{\psi_{\beta}(\beta_1 q_2)}-\frac{\psi_{\beta}(\beta_2 q_1)}{\psi_{\beta}(\beta_1 q_1)}\Big)\Big(-\frac{\psi(\beta_1 q_3)}{\psi_{\beta}(\beta_1 q_3)}+\frac{\psi(\beta_2 q_3)}{\psi_{\beta}(\beta_1 q_3)}\Big)\Big]\nonumber\\
&=\frac{\beta_2}{2\beta_{1}^{2}}\Big[\frac{\psi_{\beta}(\beta_2 q_3)}{\psi_{\beta}(\beta_1 q_3)}\Big(\frac{\psi(\beta_1 q_1)}{\psi_{\beta}(\beta_1 q_1)}-\frac{\psi(\beta_1 q_2)}{\psi_{\beta}(\beta_1 q_2)}\Big)\nonumber\\
&+\frac{\psi_{\beta}(\beta_2 q_2)}{\psi_{\beta}(\beta_1 q_2)}\Big(\frac{\psi(\beta_1 q_1)}{\psi_{\beta}(\beta_1 q_1)}-\frac{\psi(\beta_1 q_3)}{\psi_{\beta}(\beta_1 q_3)}\Big)\nonumber\\
&+\frac{\psi_{\beta}(\beta_2 q_1)}{\psi_{\beta}(\beta_1 q_1)}\Big(\frac{\psi(\beta_1 q_3)}{\psi_{\beta}(\beta_1 q_3)}-\frac{\psi(\beta_1 q_2)}{\psi_{\beta}(\beta_1 q_2)}\Big)\nonumber\\
&+\frac{\psi_{\beta}(\beta_2 q_2)}{\psi_{\beta}(\beta_1 q_2)}\cdot\frac{\psi(\beta_2 q_1)}{\psi_{\beta}(\beta_1 q_1)}+\frac{\psi_{\beta}(\beta_2 q_2)}{\psi_{\beta}(\beta_1 q_2)}\cdot \frac{\psi(\beta_2 q_3)}{\psi_{\beta}(\beta_1 q_3)}\nonumber\\
&+\frac{\psi_{\beta}(\beta_2 q_3)}{\psi_{\beta}(\beta_1 q_3)}\cdot \frac{\psi(\beta_2 q_1)}{\psi_{\beta}(\beta_1 q_1)}-\frac{\psi_{\beta}(\beta_2 q_3)}{\psi_{\beta}(\beta_1 q_3)}\cdot \frac{\psi(\beta_2 q_2)}{\psi_{\beta}(\beta_1 q_2)}\nonumber\\
&-\frac{\psi_{\beta}(\beta_2 q_1)}{\psi_{\beta}(\beta_1 q_1)}\cdot \frac{\psi(\beta_2 q_2)}{\psi_{\beta}(\beta_1 q_2)}-\frac{\psi_{\beta}(\beta_2 q_1)}{\psi_{\beta}(\beta_1 q_1)}\cdot \frac{\psi(\beta_2 q_3)}{\psi_{\beta}(\beta_1 q_3)}\Big]<0.
\end{align}
\noindent We can see that the first three terms out of nine terms  in the latter sum of inequality \eqref{equuuu} are monotonic in terms of $q.$ Next we will express the remaining six terms in similar form. We have
\begin{align}
\frac{\beta_2}{2\beta_{1}^{2}}\Big[\frac{\psi_{\beta}(\beta_2 q_3)}{\psi_{\beta}(\beta_1 q_3)}\cdot \frac{\psi(\beta_2 q_1)}{\psi_{\beta}(\beta_1 q_1)}\cdot \frac{\psi_{\beta}(\beta_2 q_1)}{\psi_{\beta}(\beta_2 q_1)}
+\frac{\psi_{\beta}(\beta_2 q_2)}{\psi_{\beta}(\beta_1 q_2)}\cdot \frac{\psi(\beta_2 q_1)}{\psi_{\beta}(\beta_1 q_1)}\cdot  \frac{\psi_{\beta}(\beta_2 q_1)}{\psi_{\beta}(\beta_2 q_1)}\nonumber\\
- \frac{\psi_{\beta}(\beta_2 q_1)}{\psi_{\beta}(\beta_1 q_1)}\cdot \frac{\psi(\beta_2 q_2)}{\psi_{\beta}(\beta_1 q_2)}\cdot \frac{\psi_{\beta}(\beta_2 q_2)}{\psi_{\beta}(\beta_2 q_2)}
-\frac{\psi_{\beta}(\beta_2 q_3)}{\psi_{\beta}(\beta_1 q_3)}\cdot \frac{\psi(\beta_2 q_2)}{\psi_{\beta}(\beta_1 q_2)}\cdot \frac{\psi_{\beta}(\beta_2 q_2)}{\psi_{\beta}(\beta_2 q_2)}\nonumber\\
+\frac{\psi_{\beta}(\beta_2 q_2)}{\psi_{\beta}(\beta_1 q_2)}\cdot \frac{\psi(\beta_2 q_3)}{\psi_{\beta}(\beta_1 q_3)}\cdot \frac{\psi_{\beta}(\beta_2 q_3)}{\psi_{\beta}(\beta_2 q_3)}
-\frac{\psi_{\beta}(\beta_2 q_1)}{\psi_{\beta}(\beta_1 q_1)}\cdot \frac{\psi(\beta_2 q_3)}{\psi_{\beta}(\beta_1 q_3)}\cdot \frac{\psi_{\beta}(\beta_2 q_3)}{\psi_{\beta}(\beta_2 q_3)}\Big]\nonumber.
\end{align}

\noindent For brevity, denote $x_j=\frac{\psi_{\beta}(\beta_2 q_j)}{\psi_{\beta}(\beta_1 q_j)},$ and $y_j=\frac{\psi(\beta_2 q_j)}{\psi_{\beta}(\beta_2 q_j)}$, where $j \in \{1,2,3\}.$ In new notation, the above expression can be expressed as

\begin{equation}
\label{eq21}
\frac{\beta_2}{2\beta_{1}^{2}}\Big(x_3 x_1 y_1 +x_2 x_1 y_1
-x_1 x_2 y_2 -x_3 x_2 y_2 
+x_2 x_3 y_3 -x_1 x_3 y_3 \Big).
\end{equation}

\noindent Adding and subtracting the terms, $\frac{\beta_2}{2\beta_{1}^{2}} x_3 x_1 y_2$ and $\frac{\beta_2}{2\beta_{1}^{2}}  x_1 x_2 y_3$ to \eqref{eq21}, we obtain

\begin{align*}
\frac{\beta_2}{2\beta_{1}^{2}}\Big(x_3 x_1 y_1 -x_3 x_1 y_2+x_3 x_1 y_2
+x_2 x_1 y_1 -x_2 x_1 y_3
 +x_2 x_1 y_3-x_1 x_2 y_2 -x_3 x_2 y_2 
+x_3 x_2 y_3-x_1 x_3 y_3\Big)\\
=\frac{\beta_2}{2\beta_{1}^{2}}\Big(x_3 x_1 (y_1 -y_2)+x_2 x_1 (y_1 -y_3)
+(y_2 -y_3)(x_3 x_1 -x_1 x_2 -x_3 x_2)\Big).
\end{align*}

\noindent Using the definition of $x_j, $ and $y_j,$ the above sum is same as 
\begin{align}
&\frac{\beta_2}{2\beta_{1}^{2}}\Big[\frac{\psi_{\beta}(\beta_2 q_3)}{\psi_{\beta}(\beta_1 q_3)}\cdot \frac{\psi_{\beta}(\beta_2 q_1)}{\psi_{\beta}(\beta_1 q_1)}\Big(\frac{\psi(\beta_2 q_1)}{\psi_{\beta}(\beta_2 q_1)}-\frac{\psi(\beta_2 q_2)}{\psi_{\beta}(\beta_2 q_2)}\Big)
+\frac{\psi_{\beta}(\beta_2 q_2)}{\psi_{\beta}(\beta_1 q_2)}\cdot \frac{\psi_{\beta}(\beta_2 q_1)}{\psi_{\beta}(\beta_1 q_1)}\Big(\frac{\psi(\beta_2 q_1)}{\psi_{\beta}(\beta_2 q_1)}-\frac{\psi(\beta_2 q_3)}{\psi_{\beta}(\beta_2 q_3)}\Big)\nonumber\\
&+\Big(\frac{\psi(\beta_2 q_2)}{\psi_{\beta}(\beta_2 q_2)}-\frac{\psi(\beta_2 q_3)}{\psi_{\beta}(\beta_2 q_3)}\Big)\Big(\frac{\psi_{\beta}(\beta_2 q_3)}{\psi_{\beta}(\beta_1 q_3)}\cdot \frac{\psi_{\beta}(\beta_2 q_1)}{\psi_{\beta}(\beta_1 q_1)}-\frac{\psi_{\beta}(\beta_2 q_1)}{\psi_{\beta}(\beta_1 q_1)}\cdot\frac{\psi_{\beta}(\beta_2 q_2)}{\psi_{\beta}(\beta_1 q_2)}-\frac{\psi_{\beta}(\beta_2 q_3)}{\psi_{\beta}(\beta_1 q_3)}\cdot \frac{\psi_{\beta}(\beta_2 q_2)}{\psi_{\beta}(\beta_1 q_2)}\Big)\Big].\nonumber
 \end{align}
 \noindent  Combining the above expression with the first three monotonic terms of inequality \eqref{equuuu}, we obtain 
 \begin{align*}
 &\frac{\beta_2}{2\beta_{1}^{2}}\Big[\frac{\psi_{\beta}(\beta_2 q_3)}{\psi_{\beta}(\beta_1 q_3)}\Big(\frac{\psi(\beta_1 q_1)}{\psi_{\beta}(\beta_1 q_1)}-\frac{\psi(\beta_1 q_2)}{\psi_{\beta}(\beta_1 q_2)}\Big)
 +\frac{\psi_{\beta}(\beta_2 q_3)}{\psi_{\beta}(\beta_1 q_3)}\cdot \frac{\psi_{\beta}(\beta_2 q_1)}{\psi_{\beta}(\beta_1 q_1)}\Big(\frac{\psi(\beta_2 q_1)}{\psi_{\beta}(\beta_2 q_1)}-\frac{\psi(\beta_2 q_2)}{\psi_{\beta}(\beta_2 q_2)}\Big)\\
 &+\frac{\psi_{\beta}(\beta_2 q_2)}{\psi_{\beta}(\beta_1 q_2)}\Big(\frac{\psi(\beta_1 q_1)}{\psi_{\beta}(\beta_1 q_1)}-\frac{\psi(\beta_1 q_3)}{\psi_{\beta}(\beta_1 q_3)}\Big)
 +\frac{\psi_{\beta}(\beta_2 q_2)}{\psi_{\beta}(\beta_1 q_2)}\cdot \frac{\psi_{\beta}(\beta_2 q_1)}{\psi_{\beta}(\beta_1 q_1)}\Big(\frac{\psi(\beta_2 q_1)}{\psi_{\beta}(\beta_2 q_1)}-\frac{\psi(\beta_2 q_3)}{\psi_{\beta}(\beta_2 q_3)}\Big)\\
 &+\frac{\psi_{\beta}(\beta_2 q_1)}{\psi_{\beta}(\beta_1 q_1)}\Big(\frac{\psi(\beta_1 q_3)}{\psi_{\beta}(\beta_1 q_3)}-\frac{\psi(\beta_1 q_2)}{\psi_{\beta}(\beta_1 q_2)}\Big)
 +\Big(\frac{\psi(\beta_2 q_2)}{\psi_{\beta}(\beta_2 q_2)}-\frac{\psi(\beta_2 q_3)}{\psi_{\beta}(\beta_2 q_3)}\Big)\Big(\frac{\psi_{\beta}(\beta_2 q_3)}{\psi_{\beta}(\beta_1 q_3)}\cdot \frac{\psi_{\beta}(\beta_2 q_1)}{\psi_{\beta}(\beta_1 q_1)}-\frac{\psi_{\beta}(\beta_2 q_1)}{\psi_{\beta}(\beta_1 q_1)}\cdot\frac{\psi_{\beta}(\beta_2 q_2)}{\psi_{\beta}(\beta_1 q_2)}\\
 &-\frac{\psi_{\beta}(\beta_2 q_3)}{\psi_{\beta}(\beta_1 q_3)}\cdot \frac{\psi_{\beta}(\beta_2 q_2)}{\psi_{\beta}(\beta_1 q_2)}\Big)\Big]<0.
 \end{align*}

\end{proof}
\end{lem}

\noindent  Denote $x_j=\frac{\psi_{\beta}(\beta_2 q_j)}{\psi_{\beta}(\beta_1 q_j)}, y_j=\frac{\psi(\beta_2 q_j)}{\psi_{\beta}(\beta_2 q_j)},$ and $z_j=\frac{\psi(\beta_1 q_j)}{\psi_{\beta}(\beta_1 q_j)}.$ In the new notation, the above inequality can be expressed as 

\begin{equation}
\label{eq22}
\frac{\beta_2}{2\beta_1 ^{2}}\Big[x_3(z_1-z_2)+x_3 x_1(y_1-y_2)+x_2(z_1-z_3)+x_2 x_1(y_1 -y_3)+x_1 (z_3-z_2)+(y_3-y_2)(x_1 x_2 +x_3 x_2-x_3 x_1)\Big]<0
\end{equation}

\indent Next, we will check whether there exist positive numbers $\beta_1,\beta_2,q_1,q_2,$ and $q_3$ such that inequality \eqref{eq22} is satisfied. We will name this problem as Three point problem.  This is the subject of discussion in the following section.

\subsection{Solution To Three Point Problem}

\indent We will present the solution to three point problem using  two different  cases. In Case I, we will assume that $q_3<\max(q_1,q_2).$ Without loss of generality we will assume that $q_1>q_2.$ The other situation, where $q_2>q_1$ can be analyzed similarly. We will show that if $q_1=\max(q_j)_{j=1}^{3},$ then \eqref{eq22} is not satisfied. To this end,  we will present some auxiliary results.

\begin{lem}
\label{lem9}
Suppose that $ q_1=\max(q_j)_{j=1}^{3}.$ Then, $\beta_1 > \beta_2$, where $g_{1}'(\beta_1)=g_{2}'(\beta_2)=0.$

\indent \begin{proof} We saw earlier that 
 $\Delta=\psi'(\beta_1 q_1)\psi'(\beta_2 q_2)-\psi'(\beta_2 q_1)\psi'(\beta_1 q_2)>0$. This implies that
 \begin{equation}
 \label{eq23}
\frac{\psi'(\beta_1 q_1)}{\psi'(\beta_1 q_2)}>\frac{\psi'(\beta_2 q_1)}{\psi'(\beta_2 q_2)}.
 \end{equation}
 
\indent Suppose by contradiction,  $\beta_1 <\beta_2.$ Moreover, assumption 2 gives $\frac{\psi'(\beta p)}{\psi'(\beta q)}, p>q$ is an increasing function of $\beta$. Using  assumption 2, and the assumption that $\beta_2 >\beta_1,$ we get   $\displaystyle \frac{\psi'(\beta_1 q_1)}{\psi'(\beta_1 q_2)}<\frac{\psi'(\beta_2 q_1)}{\psi'(\beta_2 q_2)}.$ Here $p=q_1,$ and $q=q_2.$ This contradicts  \eqref{eq23}. Hence, $\beta_1>\beta_2.$ 
\end{proof}
\end{lem}
\begin{lem}
\label{lem10}
Suppose that $q_1>q_3>q_2.$ Then the following inequalities are satisfied.
\begin{enumerate}[\rm(i)]
 \item{ $z_1 >z_3>z_2.$}
 \item{$y_1 >y_3>y_2.$}
 \item{$x_2>x_3>x_1>0.$}
 \end{enumerate}
  \begin{proof} 
  
 Notice,  $x_j=\frac{\psi_{\beta}(\beta_2 q_j)}{\psi_{\beta}(\beta_1 q_j)}, y_j=\frac{\psi(\beta_2 q_j)}{\psi_{\beta}(\beta_2 q_j)},$ and $z_j=\frac{\psi(\beta_1 q_j)}{\psi_{\beta}(\beta_1 q_j)}.$\\

 \indent \rm(i) We will show that $z_1 >z_2.$ Remaining inequalities can be shown similarly.\\\indent We need to show that $\frac{\psi(\beta_1 q_1)}{\psi_{\beta}(\beta_1 q_1)}> \frac{\psi(\beta _1 q_2)}{\psi_{\beta}(\beta_1 q_2)}$.
 From assumption 4, we know that $\frac{d}{d\beta}\Big(\frac{\psi(\beta p)}{\psi(\beta q)}\Big)<0,$ where $p>q.$ This implies that $\psi_{\beta}(\beta p)\psi(\beta q)<\psi(\beta p) \psi_{\beta}(\beta q).$ Put $\beta=\beta_1, p=q_1, $ and $q=q_2.$ Then, we obtain  
 $\psi_{\beta}(\beta_{1} q_1)\psi(\beta_{1} q_2)<\psi(\beta_{1} q_1) \psi_{\beta}(\beta_{1} q_2),$ which in turn implies that  $\frac{\psi(\beta_1 q_1)}{\psi_{\beta}(\beta_1 q_1)}> \frac{\psi(\beta _1 q_2)}{\psi_{\beta}(\beta_1 q_2)}$.
 Similarly, we can show  other inequalities. Proof for \rm(i) completed. \\
 
\indent \rm(ii)The proof is identical to (i).\\

 \indent \rm(iii)  We will show that $x_2>x_3.$ Remaining inequalities can be shown similarly.\\ \indent   We need to show that $\frac{\psi_{\beta}(\beta_2 q_2)}{\psi_{\beta}(\beta_1 q_2)}>\frac{\psi_{\beta}(\beta_2 q_3)}{\psi_{\beta}(\beta_1 q_3)}.$  Using assumption 2 and the fact that $\beta_1 >\beta_2$, we obtain  $\frac{d}{dq}\Big(\frac{\psi_{\beta}( \beta_2 q)}{\psi_{\beta}(\beta_{1} q)}\Big)<0.$ Hence we obtain $\frac{\psi_{\beta}(\beta_2 q_2)}{\psi_{\beta}(\beta_1 q_2)}>\frac{\psi_{\beta}(\beta_2 q_3)}{\psi_{\beta}(\beta_1 q_3)}.$ Moreover,  since $\psi'(x)<0,$ we get  $x_j>0$ for all $j \in\{1,2,3\}.$ This completes the proof for \rm(iii).

 \end{proof}
 \end{lem}
 \begin{lem}
 \label{lem11new}
 Suppose $q_1>q_2>q_3.$ Then the following inequalities hold true. 
 \begin{enumerate}[\rm(i)]
 \item{$x_3(y_2-y_3)<z_2-z_3.$}
 \item{$0<x_1<x_2<x_3.$}
 \item{$z_1>z_2>z_3$.}
 \item{$y_1>y_2>y_3.$}
 \end{enumerate}
 \indent \begin{proof}

 \indent Using assumption (2), we get $\frac{d}{d \beta}\Big(\frac{\psi'(\beta p)}{\psi'(\beta q)}\Big)>0, p>q.$ Recall that $\psi'(\beta q)=\frac{d}{d s}\psi(s)\vert_{s= \beta q}.$  This implies that $\frac{q}{p}\frac{d}{d\beta}\Big(\frac{\psi'(\beta p)}{\psi'(\beta q)}\Big)>0$, which, in turn implies $ \frac{d}{d\beta}\Big(\frac{\psi_{\beta}(\beta p)}{\psi_{\beta}(\beta q)}\Big)>0$. Hence $\frac{d}{d\beta}\Big(\frac{\psi_{\beta}(\beta q)}{\psi_{\beta}(\beta p)}\Big)<0$, where $p>q.$ Here $p=q_2,$ and $ q=q_3.$ Since $\psi(\cdot)$ is a non-negative function, we can rewrite the inequality $\frac{d}{d\beta}\Big(\frac{\psi_{\beta}(\beta q)}{\psi_{\beta}(\beta p)}\Big)<0$ as  
$$\psi_{\beta}(\beta q_3)-\psi_{\beta}(\beta q_2)\frac{\psi_{\beta}(\beta q_3)}{\psi_{\beta}(\beta q_2)} -\psi(\beta q_2)\frac{d}{d\beta}\Big(\frac{\psi_{\beta}(\beta q_3)}{\psi_{\beta}(\beta q_2)}\Big)>0$$ 
\indent The left hand side of the above inequality is the derivative of $\psi(\beta q_3)-\psi_{\beta}(\beta q_3)\frac{\psi(\beta q_2)}{\psi_{\beta}(\beta q_2)}.$ Since, $\beta_1>\beta_2,$ we obtain following sequence of inequalities

\begin{align}
&\psi(\beta_1 q_3)-\psi_{\beta}(\beta_1 q_3)\frac{\psi(\beta_1 q_2)}{\psi_{\beta}(\beta_1 q_2)}>\psi(\beta_2 q_3)-\psi_{\beta}(\beta_2 q_3)\frac{\psi(\beta_2 q_2)}{\psi_{\beta}(\beta_2 q_2)},\nonumber\\
&\psi_{\beta}(\beta_2 q_3)\frac{\psi(\beta_2 q_2)}{\psi_{\beta}(\beta_2 q_2)}-\psi(\beta_2 q_3)>\psi_{\beta}(\beta_1 q_3)\frac{\psi(\beta_1 q_2)}{\psi_{\beta}(\beta_1 q_2)}-\psi(\beta_1 q_3),\nonumber\\
&\frac{\psi_{\beta}(\beta_2 q_3)}{\psi_{\beta}(\beta_1 q_3)}\Big(\frac{\psi(\beta_2 q_2)}{\psi_{\beta}(\beta_2 q_2)}-\frac{\psi(\beta_2 q_3)}{\psi_{\beta}(\beta_2 q_3)}\Big)<\frac{\psi(\beta_1 q_2)}{\psi_{\beta}(\beta_1 q_2)}-\frac{\psi(\beta_1 q_3)}{\psi_{\beta}(\beta_1 q_3)}.\nonumber
\end{align}

\noindent This implies that 
$x_3(y_2 -y_3)<z_2 -z_3.$
 The proofs for remaining inequalities are similar to lemma \ref{lem10}. 
 \end{proof}
 \end{lem}

 \noindent Next, we will use the above auxiliary results to show that inequality  \eqref{eq22} is not satisfied. Notice that  $x_j=\frac{\psi_{\beta}(\beta_2 q_j)}{\psi_{\beta}(\beta_1 q_j)}, y_j=\frac{\psi(\beta_2 q_j)}{\psi_{\beta}(\beta_2 q_j)},$ and $z_j=\frac{\psi(\beta_1 q_j)}{\psi_{\beta}(\beta_1 q_j)}.$\\

 \begin{lem}
 \label{lem12}Suppose $q_1=\max(q_j)_{j=1}^{3}.$ Then inequality \eqref{eq22} does not hold.  
  \begin{proof}
 \indent \rm(i) Suppose that $q_1 >q_3>q_2.$ We need to show that 
 $$x_3(z_1-z_2)+x_3 x_1(y_1-y_2)+x_2(z_1-z_3)+x_2 x_1(y_1 -y_3)+x_1 (z_3-z_2)+(y_3-y_2)(x_1 x_2 +x_3 x_2-x_3 x_1)$$
 is positive. From (i) and (ii) in lemma \ref{lem10}, we obtain
 $$x_3(z_1-z_2)+x_3 x_1(y_1-y_2)+x_2(z_1-z_3)+x_2 x_1(y_1 -y_3)+x_1 (z_3-z_2)>0$$
 It remains to show that $(y_3 -y_2)(x_1 x_2+x_2 x_3-x_3 x_1)>0.$ From (ii) and (iii) in lemma \ref{lem10}, we obtain $y_3>y_2$ and $x_1(x_2-x_3)>0$ respectively. Hence 
   $$x_3(z_1-z_2)+x_3 x_1(y_1-y_2)+x_2(z_1-z_3)+x_2 x_1(y_1 -y_3)+x_1 (z_3-z_2)+(y_3-y_2)(x_1 x_2 +x_3 x_2-x_3 x_1)>0$$

 \indent \rm(ii)  Next suppose that $q_1>q_2>q_3.$ We need to show that 
\begin{equation}
\label{eqwwww}
x_3(z_1-z_2)+x_3 x_1(y_1-y_2)+x_2(z_1-z_3)+x_2 x_1(y_1 -y_3)+x_1 (z_3-z_2)+(y_3-y_2)(x_1 x_2 +x_3 x_2-x_3 x_1)
\end{equation}
is positive.  \\

\noindent Using parts (ii) and (iii) of lemma \ref{lem11new}, we get $x_2(z_1-z_3)>0$ and $x_1(z_3-z_2)<0.$ In addition, we also obtain that
\begin{align*}
&x_2(z_1-z_3)+x_1(z_3-z_2)=x_2(z_1-z_2+z_2-z_3)+x_1(z_3-z_2)\nonumber\\
&=x_2(z_1 -z_2)+(x_1 -x_2)(z_3 -z_2)>0.
\end{align*}
\noindent We can rewrite expression \eqref{eqwwww} as 
\begin{align*} 
&x_3(z_1-z_2)+x_3 x_1(y_1-y_2)+x_2 x_1(y_1 -y_3)\nonumber\\
&+x_2 (z_1-z_2)+(x_1-x_2)(z_3-z_2)+(y_3-y_2)(x_1 x_2 +x_3 x_2-x_3 x_1)\nonumber\\
&=x_3(z_1-z_2)+x_3 x_1(y_1-y_2)+x_1 x_2(y_1- y_2)\nonumber\\
&+x_1 x_2(y_2-y_3)+x_2 (z_1-z_2)+(x_1-x_2)(z_3-z_2)+(y_3-y_2)(x_1 x_2 +x_3 x_2-x_3 x_1)\nonumber\\
&=(z_1-z_2)(x_2+x_3)+(y_1 -y_2)(x_1 x_3 +x_1 x_2)\nonumber\\
&+(x_1 -x_2)\Big((z_3 -z_2)+x_3(y_2-y_3)\Big).
\end{align*}

\indent Using parts (ii), (iii), and (iv) of lemma \ref{lem11new} we deduce that $(z_1-z_2)(x_2+x_3)+(y_1 -y_2)(x_1 x_3 +x_1 x_2)>0.$ Next, it suffices to show that $\Big((z_3 -z_2)+x_3(y_2-y_3)\Big)<0.$ This is obvious from part (i) of lemma \ref{lem11new}.
\noindent Hence expression \eqref{eqwwww} is positive. This completes the proof for part (ii).

\indent Summarizing the results from parts (i), and (ii) we obtain that  when $q_3<\max(q_1,q_2),$ the expression  $x_3(z_1-z_2)+x_3 x_1(y_1-y_2)+x_2(z_1-z_3)+x_2 x_1(y_1 -y_3)+x_1 (z_3-z_2)+(y_3-y_2)(x_1 x_2 +x_3 x_2-x_3 x_1)$ is positive.

\end{proof}
\end{lem}

\indent Next we will consider Case II. Assume that $q_3 >max(q_1,q_2).$ Without loss of generality assume that $q_1>q_2.$ The other case can be analyzed similarly. We will show that
\begin{equation}
\label{eq26}
x_3(z_1-z_2)+ x_1 x_3(y_1-y_2)+x_2(z_1-z_3)
+x_2 x_1( y_1 -y_3)+x_1 (z_3 -z_2) 
+ (y_3 -y_2)(x_1 x_2 +x_2 x_3 -x_1 x_3) >0.
\end{equation}
Before going over the  proof, we will present some auxiliary results.  Notice that  $x_j=\frac{\psi_{\beta}(\beta_2 q_j)}{\psi_{\beta}(\beta_1 q_j)}, y_j=\frac{\psi(\beta_2 q_j)}{\psi_{\beta}(\beta_2 q_j)},$ and $z_j=\frac{\psi(\beta_1 q_j)}{\psi_{\beta}(\beta_1 q_j)}.$
\begin{lem}
\label{lem13}
Suppose $q_3>q_1>q_2.$ Then the following inequalities hold true.
\begin{enumerate}[\rm(i)]
\item{$ x_2 >x_1> x_3.$}
\item{$z_3>z_1 >z_2.$}
\item{$y_3 >y_1>y_2.$}
\item{$\frac{z_1}{x_1}>\frac{z_2}{x_2}$}
\end{enumerate}
\indent \begin{proof}
\indent\rm(i)  From assumption 2, we know that   $\frac{d}{d\beta}\Big(\frac{\psi'(\beta p)}{\psi'(\beta q)}\Big)>0, $ where $p>q.$ This  condition can be restated as $\frac{d}{dq}\Big(\frac{\psi'(\beta_2 q)}{\psi'(\beta_1 q)}\Big)<0$, where $\beta_2 <\beta_1$(lemma \eqref{lem9}). This implies that $\frac{\psi_{\beta}(\beta_2 q_2)}{\psi_{\beta}(\beta_1 q_2)}>\frac{\psi_{\beta}(\beta_2 q_1)}{\psi _{\beta}(\beta_1 q_1)}$, and we obtain $x_2 >x_1.$ Here $\psi_{\beta}(\beta p):=\frac{d}{d \beta}(\psi(\beta p)).$ \\
\noindent Similarly, we can show that $x_1 >x_3.$  Hence we get $x_2 >x_1> x_3>0.$\\

 \item \indent \rm(ii) Using assumption 4, we know that
$\frac{d}{d \beta}\Big(\frac{\psi(\beta p)}{\psi(\beta q)} \Big)<0 $ for all $p>q.$  This implies that $\psi_{\beta}(\beta q)\psi(\beta p)>\psi_{\beta}(\beta p)\psi(\beta q).$   If we put $\beta=\beta_1, p=q_1, $ and $q=q_2,$ we obtain $\frac{\psi(\beta_1 q_1)}{\psi_{\beta}(\beta_1 q_1)}>\frac{\psi(\beta _1 q_2)}{\psi_{\beta}(\beta_1 q_2)}.$ This implies $z_1 >z_2.$ Similarly, we get $z_3 >z_1.$ Hence $z_3 >z_1 >z_2.$\\

\indent \rm(iii)The proof is similar to \rm(ii).\\

\indent\rm(iv)  We need to show that $\frac{\psi(\beta_1 q_1)}{\psi_{\beta}(\beta_2 q_1)}>\frac{\psi(\beta_1 q_2)}{\psi_{\beta}(\beta_2 q_2)}.$ This is equivalent to showing that $\frac{\psi(\beta_1 q_1)}{\psi(\beta_1 q_2)}<\frac{\psi_{\beta}(\beta_2 q_1)}{\psi_{\beta}(\beta_2 q_2)}.$(since $\psi_{\beta}(\cdot)<0$)\\
\indent Using assumption (2), we obtain $\frac{d}{d\beta}\Big(\frac{\psi_{\beta}(\beta p)}{\psi_{\beta}(\beta q)}\Big)>0.$ This implies $\frac{\psi_{\beta}(\beta_2 q_1)}{\psi_{\beta}(\beta_2 q_2)}>\lim_{\beta_2 \rightarrow 0^{+}}\frac{\psi_{\beta}(\beta_2 q_1)}{\psi_{\beta}(\beta_2 q_2)}.$ Now we will show that $\lim_{\beta_2 \rightarrow 0^{+}}\frac{\psi_{\beta}(\beta_2 q_1)}{\psi_{\beta}(\beta_2 q_2)}$ exists.\\
\indent Since $\frac{d}{d\beta}\Big(\frac{\psi_{\beta}(\beta p)}{\psi_{\beta}(\beta q)}\Big)>0,$ we obtain that $\Big(\frac{\psi_{\beta}(\beta p)}{\psi_{\beta}(\beta q)}\Big)$ is a decreasing function of $\beta$, as $\beta \rightarrow 0^{+}.$ Hence $\{\frac{\psi_{\beta}(\beta p)}{\psi_{\beta}(\beta q)}\}$ has an upper bound for small but positive values of $\beta.$ Moreover since $\psi_{\beta}(\cdot)$ is negative, $\{\frac{\psi_{\beta}(\beta p)}{\psi_{\beta}(\beta q)}\}$ has a lower bound too. We obtain a monotonic and bounded sequence, hence $\lim_{\beta_2 \rightarrow 0^{+}}\frac{\psi_{\beta}(\beta_2 q_1)}{\psi_{\beta}(\beta_2 q_2)}$ exists.\\
\indent  Using L'hospital's rule we get $\lim_{\beta_2 \rightarrow 0^{+}}\frac{\psi_{\beta}(\beta_2 q_1)}{\psi_{\beta}(\beta_2 q_2)}=\lim_{\beta_2 \rightarrow 0^{+}}\frac{\psi(\beta_2 q_1)}{\psi(\beta_2 q_2)}.$  Notice that $\psi(y)\rightarrow \infty, $ as $y \rightarrow 0^{+}.$ In addition, assumption (4) gives $\frac{d}{d\beta}\Big(\frac{\psi(\beta p)}{\psi(\beta q)}\Big) <0, p>q.$ This implies that $\lim_{\beta_2 \rightarrow 0^{+}}\frac{\psi(\beta_2 q_1)}{\psi(\beta_2 q_2)}>\frac{\psi(\beta_1 q_1)}{\psi(\beta_1 q_2)}.$  Hence we obtain the following chain of  inequalities 

$$\frac{\psi_{\beta}(\beta_2 q_1)}{\psi_{\beta}(\beta_2 q_2)}>\lim_{\beta_2 \rightarrow 0^{+}}\frac{\psi_{\beta}(\beta_2 q_1)}{\psi_{\beta}(\beta_2 q_2)}\overset{LH}{=}\lim_{\beta_2 \rightarrow 0^{+}}\frac{\psi(\beta_2 q_1)}{\psi(\beta_2 q_2)}>\frac{\psi(\beta_1 q_1)}{\psi(\beta_1 q_2)}$$

\noindent Therefore we get $\frac{\psi(\beta_1 q_1)}{\psi(\beta_1 q_2)}<\frac{\psi_{\beta}(\beta_2 q_1)}{\psi_{\beta}(\beta_2 q_2)}.$ Proof for \rm(iv) completed.

\end{proof}
\end{lem}

\begin{lem}
\label{lem14}
Suppose that $q_{3}>q_1>q_2.$ Then inequality \eqref{eq26} holds true. 
\indent \begin{proof} From 
 lemma \ref{lem9}, we obtain that  $\beta_1>\beta_2.$ The left hand side of \eqref{eq26} can be rewritten as 
\begin{align}
&x_3(z_1-z_2)+ x_1 x_3(y_1-y_2)+x_2 z_1 -x_2 z_3\nonumber\\
&+x_1(z_3-z_2)+x_1 x_2(y_1 -y_2 +y_2 -y_3)\nonumber\\
&+(y_3 -y_2)x_1 x_2 +(y_3 -y_2)(x_2 x_3-x_1 x_3)\nonumber\\
&=x_3(z_1-z_2)+ x_1 x_3(y_1-y_2)+z_3(x_1 -x_2)\nonumber\\
&+x_2 z_1- x_1 z_2 +x_1 x_2(y_1 -y_2)+(y_3 -y_2)x_3 (x_2 -x_1).\nonumber
\end{align}

\noindent Using parts (i), (ii), and (iii), of lemma \ref{lem13},  we get $x_3(z_1 -z_2)>0, x_1 x_3(y_1-y_2)>0,z_3(x_1 -x_2)>0,$ and $(y_3 -y_2)x_3 (x_2 -x_1)>0.$\\
\noindent Hence $x_3(z_1-z_2)+ x_1 x_3(y_1-y_2)+z_3(x_1 -x_2)
+x_2 z_1- x_1 z_2 +x_1 x_2(y_1 -y_2)+(y_3 -y_2)x_3 (x_2 -x_1)
>x_2 z_1 -x_1 z_2+(x_1 x_2)(y_1 -y_2).$
\indent Next we will show that $x_2 z_1 -x_1 z_2+(x_1 x_2)(y_1 -y_2)>0.$ Since $x_1 x_2 (y_1 -y_2)>0$(using(iii)), it suffices to show that $x_2 z_1 -x_1 z_2>0.$ This is equivalent to showing that $\frac{z_1}{x_1}>\frac{z_2}{x_2}.$  But we saw in (iv) that $\frac{z_1}{x_1}-\frac{z_2}{x_2}>0.$ Hence we obtain 
\begin{align}
&x_3(z_1-z_2)+ x_1 x_3(y_1-y_2)+x_2 z_1 -x_2 z_3\nonumber\\
&+x_1(z_3-z_2)+x_1 x_2(y_1 -y_2 +y_2 -y_3)\nonumber\\
&+(y_3 -y_2)x_1 x_2 +(y_3 -y_2)(x_2 x_3-x_1 x_3)>0.\nonumber
\end{align}
This completes the proof for the case when $q_3>q_1>q_2.$
\end{proof}
\end{lem}

\indent This completes the solution to Three Point problem. Now we have all the required tools to present the main result of this section. We will show that the function $f(x)=\sum_{i=1}^{n}c_i\phi(x_i), c_i \neq 0$ for all $i$ has at most one point of local maxima in two side \octants.\\

\subsection{Number of Points of Local maxima in two Side orthants}

\indent Recall that in the first side \octant \; $f(x)$ is of the form 
$g_{1}(\beta)=\sum_{j=2}^{n}l_j \psi(\beta q_j)-l_1 \psi(\beta q_1),$
and in the second side \octant \; $f(x)$ is of the form $g_{2}(\beta)=\sum_{j=1,j \neq 2}^{n}l_j \psi(\beta q_j)-l_{2}\psi(\beta q_2).$ It can be easily checked that 
$g_{1}'(\beta)=\sum_{j=2}^{n}l_j q_j \psi'(\beta q_j)-l_1 q_1 \psi'(\beta q_1),$
and $g_{2}'(\beta)=\sum_{j=1,j \neq 2}^{n}l_j q_j \psi'(\beta q_j)-l_2 q_2 \psi'(\beta q_2).
$\\
\indent We will first show that for small values of $\beta,$ either $g_{1}(\beta)$ or $g_{2}(\beta)$ is a non-increasing function. Since $\psi'(x)<0$, we obtain that $\sign(g_{1}'(\beta))=-\sign\Big(\sum_{j=2}^{n} l_j q_j \frac{\psi'(\beta q_j)}{ \psi'(\beta q_1)}-l_1 q_1\Big).$ Similarly for $g_{2}'(\beta)$.\\
\indent Suppose $\lim_{\beta \rightarrow 0}g_{1}'(\beta)>0,$ and $\lim_{\beta \rightarrow 0}g_{2}'(\beta)>0.$ This implies that $\Big(\sum_{j=2}^{n} l_j q_j \frac{\psi'(\beta q_j)}{ \psi'(\beta q_1)}-l_1 q_1\Big)<0$ for values of beta close to zero. Since $l_j,$ and $q_j$ are positive for all $j$, we get 
$l_1 q_1 >l_2 q_2 \frac{\psi'(\beta q_2)}{\psi'(\beta q_1)}
$. Similarly, $ \lim_{\beta \rightarrow 0}g_{2}'(\beta)>0$ implies that $l_2 q_2 >l_1 q_1 \frac{\psi'(\beta q_1)}{\psi'(\beta q_2)}
$. From above inequalities we obtain 
$\frac{l_1 q_1}{l_2 q_2}>\frac{\psi'(\beta q_2)}{\psi'(\beta q_1)}$
and, $\frac{l_1 q_1}{l_2 q_2}<\frac{\psi'(\beta q_2)}{\psi'(\beta q_1)}$. Hence contradiction. \\

\begin{prop}
\label{prop10}
The function $f(x)=\sum_{i=1}^{n}c_i \phi(x_i)$ does not have points of local maxima in both side \octants.
\begin{proof} 
We showed earlier that for small values of $\beta,$ either $g_{1}(\beta)$ or $g_{2}(\beta)$ is non-increasing function. In this proof, we will only consider the situation, when $\lim_{\beta \rightarrow 0}g_{1}'(\beta)<0.$ Other case can be analyzed similarly.

\indent \underline{Case I}: Suppose $q_1 >q_j$ for all $j \in \{2 \ldots n\}.$  Then  $g_{1}'(\beta ) \rightarrow \infty, $ and $g_{2}'(\beta)\rightarrow -\infty,$ as $\beta \rightarrow \beta_{max}.$ From remark \ref{rem2}, we obtain that $g_{1}(\cdot)$ has at most one stationary point, and $g_{2}(\cdot)$ has at most two stationary points. These observations will be used throughout the proof of this case.  We will consider the following sub cases.

\indent \underline{Sub Case \rm(i)}: Suppose $g_{2}'(\beta)<0 $ for all $\beta$, and $\lim_{\beta \rightarrow 0}g_{1}'(\beta)<0.$  Using the fact $g_{1}'(\beta) \rightarrow \infty, $ as $\beta \rightarrow \beta_{max}$, we deduce that $g_{1}(\beta)$ changes monotonicity exactly once. So for a given $ b \in \mathbb{R},$ there exists at most one  $\beta_1 \in [0,\beta_{max}]$ such that $g_{1}(\beta_1)=b, $ and $g_{1}'(\beta_1)>0.$ Using corollary \ref{cor2}, we conclude that, $f(x)$ has at most one point of local maximum in the side \octant \; with first component negative.\\
\indent In addition, $g_{2}'(\beta)<0$ for  all $\beta$ implies that  the function $f(x)$ does not have any point of local maximum in the \octant \; with second component negative ( corollary \ref{cor2}). Hence $f(x)$ has at most one point of local maximum in  the side \octants. This completes the proof for \rm(i).

\indent \underline{Sub Case \rm(ii)}: Suppose $\lim_{\beta \rightarrow 0}g_{1}'(\beta)<0,$ and $\lim_{\beta \rightarrow 0}g_{2}'(\beta)<0.$ We know that the function $g_{1}(\beta)$ has at most one critical point in the interval $(0,\beta_{max}].$ Together with the fact that  $g_{1}'(\beta) \rightarrow \infty, $ as $\beta \rightarrow \beta_{max}$, we deduce that $g_{1}(\beta)$ changes monotonicity once. Denote critical point of $g_{1}(\beta)$ as $\beta_1.$ It can be easily checked that  $g_{1}(\beta)$ is increasing on the interval $(\beta_1, \beta_{max}]$, and decreasing otherwise.\\
\indent Meanwhile, since  $g_{2}'(\beta) \rightarrow -\infty ,$ as $\beta \rightarrow \beta_{max}$, and $\lim_{\beta \rightarrow 0} g_{2}'(\beta)<0$, we conclude that  $g_{2}(\beta)$ is either monotonically decreasing or changes monotonicity.  If $g_{2}(\beta)$ is monotonically decreasing, then we get the situation similar to Sub Case \rm(i). Here, we will consider the situation when $g_{2}(\beta)$ changes monotonicity. It is easy to see that function $g_{2}(\beta)$ has two critical points, given by $\beta'', $ and $\beta_2.$ It can be easily checked that  $g_{2}'(\beta)>0$ for all $\beta \in (\beta'',\beta_2)$, and $g_{2}'(\beta)<0$ elsewhere. \\
\indent Suppose by contradiction, $f(x)$ has two points of local maxima, one in each \octant. Hence there exist $\beta' \neq \beta'''$ such that for a given $b \in \mathbb{R}$,  we have $g_{1}(\beta')=g_{2}(\beta''')=b,g_{1}(\beta')>0,$ and $g_{2}(\beta''')>0.$ Here $\beta' \in (\beta_1, \beta_{max}),$ and $\beta''' \in (\beta'', \beta_2).$\\
\indent Since $g_{1}(\cdot)$ is increasing on the interval $(\beta_1 ,\beta_{max}),$ and $\beta' \in (\beta_1, \beta_{max}),$ we get $g_{1}(\beta_1)<g_{1}(\beta')=b$. Similarly, we obtain   $g_{2}(\beta_2)>b.$ This implies that $g_{1}(\beta_1)-g_{2}(\beta_2)<0.$ This leads us to the Three point problem. See Lemma \eqref{lem8}.Hence we get 
$$x_3(z_1-z_2)+ x_1 x_3(y_1-y_2)+x_2(z_1-z_3)
+x_2 x_1( y_1 -y_3)+x_1 (z_3 -z_2) 
+ (y_3 -y_2)((x_1 x_2 +x_2 x_3 -x_1 x_3) <0.
$$

\indent But in Lemma \ref{lem12}, we showed that if $q_{1}=\max(q_j)_{j=1}^{3},$ then 
$$x_3(z_1-z_2)+ x_1 x_3(y_1-y_2)+x_2(z_1-z_3)
+x_2 x_1( y_1 -y_3)+x_1 (z_3 -z_2) 
+ (y_3 -y_2)((x_1 x_2 +x_2 x_3 -x_1 x_3) >0.
$$

Hence contradiction  to our assumption. Therefore, $f(x)$ cannot have points of local maxima in both side \octants. Proof for \rm(ii) completed.

\indent \underline{Sub Case \rm(iii)}: Suppose $\lim_{\beta \rightarrow 0}g_{1}'(\beta)<0, $ and $\lim_{\beta \rightarrow 0}g_{2}'(\beta)>0.$  We saw earlier that $g_{1}'(\beta)\rightarrow \infty, $ as $\beta \rightarrow \beta_{max}.$ Since $g_{1}'(\beta)$ has at most one root, there exists unique $\beta_1 \in (0,\beta_{max})$ such that $g_{1}'(\beta_1)=0 $. It can be easily checked that $g_{1}'(\beta)<0 $ in the interval $ (0,\beta_1)$ and $g_{1}'(\beta)>0 $ in the interval $(\beta_1 ,\beta_{max}].$ \\
\indent At the same time,  $g_{2}'(\beta)\rightarrow -\infty, $ as $\beta \rightarrow \beta_{max}.$ Since $g_{2}'(\beta)$ has at most  two roots,  there exists unique $\beta_2 \in (0, \beta_{max})$ such that $g_{2}'(\beta_2)=0.$ It is easy to see that $g_{2}'(\beta)>0 $ in the interval $ (0,\beta_2)$ and $g_{2}'(\beta)<0 $ in the interval $(\beta_2 ,\beta_{max}).$\\
\indent Assume that $f(x)$ has two points of local maxima, one in each side \octant. Using corollary \ref{cor2}, we conclude that for some $b \in\mathbb{R},$ there exists $\beta' \neq \beta^{''} $ such that $g_{1}(\beta')=g_{2}(\beta'')=b>0,g_{1}'(\beta')>0,$ and $g_{2}'(\beta'')>0.$ It can be easily checked that $\beta' \in (\beta_1 ,\beta_{max}),$ and $\beta'' \in (0,\beta_2).$ \\
\indent Notice that $g_{1}(\beta_1)<b$ and $g_{2}(\beta_2)>b.$ Hence $g_{1}(\beta_1)-g_{2}(\beta_2)<0.$ This case can be analyzed in similar manner as Sub Case \rm(ii).  We will obtain a contradiction to the fact that $g_{1}(\beta_1)-g_{2}(\beta_2)<0.$ Hence $f(x)$ has at most one point of local maximum in both side \octants. This completes the proof for \rm(iii).

\indent To recapitulate, when $q_1 >q_j$ for all $j \in \{2 \ldots n\},$ the function $f(x)$ has at most one point of local  maximum in two side \octants. This completes the proof for Case I.

\indent {Case II}: Let $q_2=\max(q_j)_{j=1}^{n}.$ The proof for this case is similar to Case I.

\indent {Case III}. Let $q_{j0}=\max(q_j)_{j=1}^{n}.$ Here $j0 \neq 1, 2.$ Without loss of generality, let $j0=3.$ Then $g_{1}'(\beta) \rightarrow -\infty,$ and  $g_{2}'(\beta) \rightarrow -\infty, $ as $\beta \rightarrow \beta_{max}$.  Moreover, using  remark \ref{rem2}, we obtain that both functions $g_{1}(\cdot)$ and $g_{2}(\cdot)$ have at most two stationary points. These observations will be used throughout the proof of this case.
 We will consider the  following sub cases.
  
\indent  \underline{Sub Case \rm(i)}: Suppose $g_{1}'(\beta) <0, $ and $g_{2}'(\beta) <0$ for all values of $\beta \in (0,\beta_{max})$.  Using corollary  \ref{cor2}, $f(x)$ does not have any point of local maximum in both side \octants.
  
  \indent \underline{Sub Case \rm(ii)}: Assume that  $\lim_{\beta \rightarrow 0}g_{1}'(\beta)<0$, and $g_{2}'(\beta)<0$ for all $\beta \in (0,\beta_{max}).$ If $g_{1}'(\beta)$ is negative  for all $\beta,$ then we obtain same situation as Sub Case \rm(i) above.\\
\indent Here, we will assume that that $g_{1}(\beta)$ changes monotonicity. We know that, $g_{1}(\beta)$ can have  at most two critical points in $(0,\beta_{max})$. Moreover, $g_{1}'(\beta) \rightarrow -\infty, $ as $\beta \rightarrow \beta_{max}.$ Hence there exist $\beta' \neq \beta''$ such that $g_{1}'(\beta')=g_{1}'(\beta'')=0.$ Notice that $g_{1}(\beta)$ is increasing  on the interval $(\beta', \beta'')$ and is decreasing elsewhere. This implies that  for given $b \in \mathbb{R}$ there exists at most one value of $\beta,$ say $\beta_3$ such that $g_{1}(\beta_3)=b, $ and $g_{1}'(\beta_3)>0.$ It is easy to see that  $\beta_3 \in (\beta',\beta'').$\\
\indent Next, $g_{2}'(\beta)<0$ for all $\beta \in (0,\beta_{max})$ implies that $f(x)$  does not have any point of local maxima in the side \octant \; with second component negative (corollary \ref{cor2}). \\
\indent From the above discussion, we conclude that $f(x)$ has at most one point of local maximum in side \octant \; with first component negative. This completes the proof for \rm(ii).
 
 \indent \underline{Sub Case \rm(iii)}: Assume that $\lim_{\beta \rightarrow 0}g_{1}'(\beta)<0,$ and $\lim_{\beta \rightarrow 0}g_{2}'(\beta)>0.$ Since $g_{2}'(\beta) \rightarrow -\infty $ as $\beta \rightarrow \beta_{max}$, we deduce that either $g_{2}(\beta)$ has one stationary point or three stationary points. But we know that $g_{2}(\beta)$ can have at most two stationary points. Therefore there exists $\beta_2 \in (0,\beta_{max})$ such that $g_{2}'(\beta_2)=0.$ It can be easily checked that $g_{2}'(\beta)>0 $ in $(0, \beta_2)$ and $g_{2}'(\beta) <0 $ in $(\beta_2,\beta_{max}).$ \\
\indent In the side \octant \; with first component negative, we can have two possibilities, either $g_{1}'(\beta)<0$ for all $\beta \in (0,\beta_{max}),$ or $g_{1}(\beta)$ changes monotonicity.\\
 \indent If   $g_{1}'(\beta)<0$ for all $\beta \in (0,\beta_{max}),$ then $f(x)$ does not have any point of local maximum in side \octant \; with first component negative ( corollary \ref{cor2}). Moreover, since $g_{2}(\beta)$ has exactly one critical point, there exists at most one value of $\beta \in (0,\beta_2)$ such that for given $b \in \mathbb{R},g_{2}(\beta)=b,$ and $ g_{2}'(\beta)>0.$ This implies that $f(x)$ has at most one point of local maximum in side \octant \; with second component negative. \\
\indent Next, suppose $g_{1}(\cdot)$ changes monotonicity. Since $g_{1}'(\beta) \rightarrow -\infty, $ as $ \beta \rightarrow \beta_{max},$ the function $g_{1}(\beta) $ will have exactly two stationary points. 
 Denote the two stationary points, as $\beta_1,$ and $\beta'',$ where  $\beta _1 < \beta''.$ Notice that 
  $g_{1}'(\beta)>0 $ in the interval $(\beta_1, \beta'')$ and is negative elsewhere. In addition, $g_{2}'(\beta)>0 $ in $(0, \beta_2)$.\\
\indent Suppose $f(x)$ has two points of local maxima, one in each side \octant. This implies there exist $ \beta_3 \neq \beta_4 $ such that for a given $b \in \mathbb{R},$ we have $g_{1}(\beta_3)=g_{2}(\beta_4)=b>0, g_{1}'(\beta_3)>0$, and $g_{2}'(\beta_4)>0.$  It is easy to see that $\beta_4 \in (0,\beta_2)$ and $\beta_3 \in (\beta_1 ,\beta'').$\\
\indent Since $g_{1}'(\beta)>0 $ for all $\beta \in (\beta_1 ,\beta'')$ and $g_{1}(\beta_3)=b,$ we get $g_{1}(\beta_1)<b.$ Similarly we obtain $g_{2}(\beta_2)>b.$ Hence we obtain $g_{1}(\beta_1)-g_{2}(\beta_2)<0.$ This leads to the Three Point problem, which in turn implies that 
 $$x_3(z_1-z_2)+ x_1 x_3(y_1-y_2)+x_2(z_1-z_3)
+x_2 x_1( y_1 -y_3)+x_1 (z_3 -z_2) 
+ (y_3 -y_2)((x_1 x_2 +x_2 x_3 -x_1 x_3) <0.
$$

\noindent But in lemma \ref{lem14}, we showed that if $q_{3}=\max(q_j)_{j=1}^{3},$ then 
$$x_3(z_1-z_2)+ x_1 x_3(y_1-y_2)+x_2(z_1-z_3)
+x_2 x_1( y_1 -y_3)+x_1 (z_3 -z_2) 
+ (y_3 -y_2)((x_1 x_2 +x_2 x_3 -x_1 x_3) >0.
$$

Hence contradiction  to our assumption. This completes the proof for Case III.  Therefore, $f(x)$ does not  have points of local maxima in both side \octants. 

\end{proof}
\end{prop}

\section{Special Case}
\noindent In previous sections, we have analyzed the behavior of function $f(x)=\sum_{i=1}^{n}c_i \phi(x_i)$ only in the open orthant (i.e. $x_j \neq 0 $ for all $j \in \{1\ldots n\}$). Moreover, we also noticed that $f(x)$ can have at most one point of local maximum in the open orthant on hyperplane P.   But it might happen that function $f(\cdot)$ has two points of local maxima, one in the open orthant and other on the boundary of the orthant (i.e. $x_j=0 $ for at least one $j \in \{1\ldots n\}$). In this section, we will show that this cannot happen.\\
\indent First, notice that the function $f(x)$ cannot have infinitely many points of local maxima on the boundary of the orthant. Suppose $f(x)$ has infinitely many points of local maxima, $\{x_i\}_{i=1}^{\infty}.$ Since $x_i$ lies on the boundary, then for all $i$, we obtain  $x_{i}^{j}=0$ for at least one $j \in \{1\ldots n\}.$ This implies that $g(\beta_k)=b$, for the sequence $\{\beta_k\}_{k=1}^{\infty}$, where $g(\beta):=\sum_{j=1}^{n}\pm l_j \psi(\beta \frac{l_j}{c_j}).$  Recall that $\psi(\cdot)=\phi'(\cdot)^{-1}.$ Using the fact that $\phi(\cdot)$ is an analytic function, we obtain that $g(\beta)\equiv b.$ This implies that $g(\cdot)$ is constant function. Hence contradiction. Therefore, the points of local maxima on the boundary of the orthant are isolated.\\

\begin{lem}
\label{lemnew}
Suppose $x_0$ is a isolated point of local maximum for function $f(\cdot)$ on hyperplane $P=\{x:l^{T}x=b\}.$ Then, for all $\varepsilon>0$ there exists $\delta_0>0,  \tilde{l}, \tilde{b} $ such that $\Vert l-\tilde{l}\Vert<\delta_0, \vert b-\tilde{b}\vert <\delta_0$, and  a point $\tilde{x_0}$ with all components nonzero such that $\tilde{x_0}$ is a point of local maximum of function $f(\cdot)$ over $\tilde{P}=\{x:\tilde{l}x=\tilde{b}\}\cap \Vert x_0 -\tilde{x_0}\Vert <\varepsilon.$    
\begin{proof} 
Since $x_0$ is a isolated point of local maximum for $f(\cdot) $ on P, there exists sufficiently small  $\varepsilon_0>0,$ such that for all $x,$ where $x \in P$ and $\Vert x-x_0\Vert <\varepsilon_0,$ we get $f(x)\leq f(x_0).$ \\
\noindent  Fix a positive number $\varepsilon < \varepsilon_0.$ Since $x_0$ is a point of local maximum, there exists $\beta \in \mathbb{R}^{+}$ such that $x_{0}^{j}=\pm \psi(\beta \frac{l_j}{c_j})$ for all $j \in \{1\ldots n\}.$ If all components of $x_0$ are nonzero, then we are done. Suppose that $x_0$ has zero components. Without loss of generality, let $x_{0}^{1}=\ldots=x_{0}^{m}=0.$ Then, using assumption (1), we obtain $\beta\frac{l_k}{c_k}=\phi'(0)$ for all $k \in \{1\ldots m\}.$\\
\indent Since $f(\cdot)$ is continuous function, there exists $\delta_0>0$ such that $f(\cdot)$ has a point of local maximum over the set $\{x:\tilde{l}x=\tilde{b}\} \cap \{x:\Vert x-x_0\Vert <\varepsilon\}$ for all vectors $\tilde{l}, $ and numbers $\tilde{b}$ such that $\Vert l-\tilde{l}\Vert <\delta_0,$ and $\vert b-\tilde{b}\vert <\delta_0.$ \\
\noindent Fix positive numbers $\delta, \delta_1,$ such that $\delta_1<\delta.$ Denote $\tilde{l}_1=l_1 -\delta_1, \tilde{l}_j=l_j -\delta, $ for $2 \leq j \leq m, \tilde{l}_j=l_j$ for $j>m.$ Denote by $\tilde{x}_0$ the point with coordinates $\tilde{x_0}^{1}=\pm \psi(\beta \frac{l_1 -\delta_1}{c_1}), \tilde{x_0}^{j}=\pm \psi(\beta \frac{l_j -\delta}{c_j})$ for $2\leq j ,m,\tilde{x_0}^{j}=x_{0}^{j}$ for $j>m.$ Denote $\tilde{b}:=\tilde{l}^{T}\tilde{x_0}.$ Assume that $\frac{\tilde{l}_1}{c_1}<\frac{\tilde{l}_j}{c_j}$ for all $j \in \{2\ldots n\}.$ \\
\indent We will choose  $\delta$ such that $\Vert l-\tilde{l}\Vert<\delta_0,$ and $\vert b-\tilde{b}\vert <\delta_0.$ Then $f(\cdot)$ has  a point of  local maximum, $\tilde{x}$ over the set $\tilde{P}=\{x:\tilde{l}^{T}x=\tilde{b}\}\cap \{ x:\Vert x-x_0\Vert <\varepsilon\}.$ There exists $\tilde{\beta}>0$ such that $\tilde{x}^{j}=\pm \psi(\tilde{\beta}\frac{\tilde{l}_j}{c_j})$ for all $j \in \{1\ldots n\}.$ If all components of $\tilde{x}$ are nonzero, then we are done.  \\
 \indent Assume at least one component of $\tilde{x}$ is zero. Due to choice of $\delta_1,$ it should be the first component $\tilde{x}^1.$ This implies that $\tilde{\beta }\frac{\tilde{l}_1}{c_1}=\phi'(0)$, and we get $\tilde{l}\tilde{x}=\tilde{b}.$ Next consider positive number $\hat{\delta}$ and a vector $\hat{l}$ with components $\hat{l}_1=\tilde{l}_1, \hat{l_j}=\tilde{l}_j -\hat{\delta}$ for $j \in \{2 \ldots n\}.$ Define the vector $\hat{x}$ with components $\hat{x}^j=\pm \psi(\tilde{\beta }\frac{\hat{l}_{j}}{c_j})$, where $j \in \{1 \ldots n \}.$ The number $\hat{\delta}$ is chosen in such way that $\Vert \hat{l}-l\Vert <\delta_0, \hat{l}^{T} \hat{x}\neq \tilde{b}$, and $\frac{\hat{l}_1}{c_1}>\frac{\hat{l}_j}{c_j}$ for $j \in \{ 2\ldots n\}.$ Consider the hyperplane $\hat{P}=\{x:\hat{l}^{T}x=\tilde{b}\}$. There exists a point of local maximum, $y$ of function over $\hat{P}\cap \Vert y-x_0\Vert <\varepsilon.$ Therefore there exists $\hat{\beta} \in \mathbb{R}^{+}$ such that $y^j=\pm \psi(\hat{\beta}\frac{\hat{l}_j}{c_j})$ for all $j=1\ldots n.$ Since $\hat{l}^{T} \hat{x}\neq \tilde{b}$, we get that $\hat{x} \notin \hat{P}.$\\
\indent If all components of $y$ are nonzero, then we are done. Assume $y$ has a zero component. Then due to choice of $\delta_1,\hat{\delta}$, we get $y^1=0, y^j\neq 0$ where$j \in \{2\ldots n\}.$ This implies that $\frac{\hat{\beta}\hat{l}_1}{c_1}=\frac{\tilde{\beta} \hat{l}_1}{c_1}=\phi'(0).$ Therefore $\hat{\beta}=\tilde{\beta}. $ But in this case $y=\hat{x},$ and therefore vector $y$ does not belong to hyperplane $\hat{P}.$ Hence we obtain contradiction. Therefore components of vector $y$ are nonzero. Thus $y$ is a point of local maximum of function $f$ over $\hat{P}$ such that $\Vert \hat{l}-l\Vert <\delta_0, \vert b-\hat{b}\vert <\delta_0,$ and $\Vert y-x_0\Vert <\varepsilon.$

\end{proof}
\end{lem}

\noindent As a consequence we obtain the following result.\\

\begin{cor}
The function $f(x)$ does not have two points of local maximum; one in open orthant and other on the boundary of orthant.
\begin{proof}
Suppose $x_0,$ and $x_1$ are two points of local maxima for function $f(\cdot).$ Without loss of generality we can assume that $x_0$ is in the open orthant i.e. $x_{0}^{j}\neq 0 $ for all $j \in \{1\ldots n\}.$  Since $f(\cdot)$ has at most one point of local maximum in the open orthant, $x_0$ is isolated point. Then using lemma \ref{lemnew} we can find another point of local maximum $\hat{x_{0}}$ such that all the components of $\hat{x_{0}}$ are nonzero. \\
\indent Using similar reasoning as above we can also find another point of local maximum $\hat{x_{1}}$ such that all the components of $\hat{x_{1}}$ are nonzero. Then we end up with two points of local maxima in the open orthant,  on the new hyperplane $\hat{P}.$ Hence contradiction. 
\end{proof}
\end{cor}

\section{Main Result}
\noindent In  Section 3, some assumptions about the function $\phi(\cdot)$  were listed.
 Those assumptions turned out to be sufficient conditions for the above results to be true. Therefore, we obtain the following result.\\
  
  \begin{theorem}
  \label{thm7} Consider the function $f(x)=\sum_{i=1}^{n}c_i \phi(x_i), $ where $c_i \neq 0.$ Suppose the function   $\phi(\cdot) $ satisfies the conditions $(1)-(5)$ stated in Section 3.
  The function $f(x)$ has at most one point of local maximum on the hyperplane, P defined by $\{x:l'x=b\}$ for some $b \in \mathbb{R}.$
\begin{proof}
The proof follows from the previous results.
\end{proof}
\end{theorem}

\noindent In the next section, we will present couple of examples of neuron transfer functions, $\phi(\cdot),$ which satisfy the conditions proposed in section 3.

\section{Examples}
 \noindent In section 3, some assumptions regarding function $\phi(\cdot)$ were presented. In the previous section, we saw that these assumptions turned out to be sufficient conditions for Theorem \ref{thm7} to be true. In this section we will check these properties for two functions, namely $\tanh(\cdot)$ , and $\arctan(\cdot).$\\

\begin{enumerate}
\item{$\phi(\cdot) \in C^2,\phi(-x)=-\phi(x),\phi'(x)>0, x\phi''(x)<0,$ for all $x\neq 0,$ and $\lim_{x \rightarrow \infty}\phi(x)<\infty.$\\

\noindent\begin{enumerate}[(i)] \item{Suppose $\phi(\cdot)=\tanh(\cdot).$ Notice that $\phi''(x)= -\frac{\arctanh (x)}{\cosh^{2}(x)}.$ It is easy to see the the graph of $\phi''(\cdot)$ lies in second and fourth quadrants. Remaining properties can be easily checked.}

\item{\noindent Now we will consider the case when $\phi(\cdot)=\arctan(\cdot).$ Notice that $\phi''(x)=-\frac{x}{(1+x^2)^{2}}.$ Remaining properties can be checked easily.}
\end{enumerate}}
\item{The function $\phi'(\cdot)$ is invertible. \\

\begin{enumerate}[(i)]\item{\noindent Suppose that $\phi(x)=\tanh (x).$ Then $\phi'(x)=\sech^{2}(x)$, where $x\geq 0.$   Denote $\sech^{2}(x)=y.$ This implies $\psi(x)=\arctanh(\sqrt{1-x})$, where $\psi(\cdot):=\phi^{-1}(\cdot)$.} \\

\item{ Next, consider $\phi(x)=\arctan(x).$ This implies $\phi'(x)=\frac{1}{1+x^2}.$ Hence $\psi(x)=\sqrt{\frac{1}{x}-1}.$}
\end{enumerate}}

\item{$\psi(x)$ is decreasing function of $x.$\\
\begin{enumerate}[(i)]
\item{
\noindent Suppose  $\psi(x)= \arctanh (\sqrt{1-x}).$ Using definition of $\arctanh(x),$ we get $ \displaystyle\arctanh (\sqrt{1-x})=\frac{1}{2}\ln \Big(\frac{1+\sqrt{1-x}}{1-\sqrt{1-x}}\Big).$ Now, we will evaluate $\psi'(x).$

\begin{align}
\psi'(x)&=\frac{1}{2x}\Big( \frac{-1+\sqrt{1-x}}{2\sqrt{1-x}}-\frac{1+\sqrt{1-x}}{2\sqrt{1-x}}\Big)\nonumber\\
&= -\frac{1}{2x\sqrt{1-x}}.\nonumber
\end{align}}
\item{

\noindent Similarly, for the case of $\psi(x)=\sqrt{\frac{1}{x}-1},$ we obtain $\psi'(x)=-\displaystyle \frac{1}{2x^{\frac{3}{2}}\sqrt{1-x}}<0.$}
\end{enumerate}}

\item{$x (\ln\vert \psi'(x)\vert )'$ is a monotonically increasing function of $x.$\\

\begin{enumerate}[(i)]
\item{
\noindent First, consider $\psi(x)=\arctanh (\sqrt{1-x}).$ Then 
\begin{align}
(\ln \vert \psi'(x)\vert)'&= \frac{1}{2x\sqrt{1-x}}\cdot -2\Big(\sqrt{1-x}-\frac{x}{2\sqrt{1-x}}\Big)\nonumber\\
&=\frac{3x-2}{2x(1-x)}.\nonumber
\end{align}

\noindent This implies that 
$\displaystyle \frac{d}{dx}(x( \ln \vert \psi'(x)\vert)')=\frac{1}{2}\Big(\frac{3-3x+3x-2}{(1-x)^{2}}\Big)=
\frac{1}{2(1-x)^{2}}>0.$

\noindent Therefore $x (\ln\vert \psi'(x)\vert )'$ is a monotonically increasing function of $x.$\\}
\item{

\noindent Next, suppose $\psi(x)=\sqrt{\frac{1}{x}-1}.$ Then $\psi'(x)=\displaystyle -\frac{1}{2x^{\frac{3}{2}}\sqrt{1-x}}.$ Now we will evaluate $(\ln \vert \psi'(x)\vert)'.$\\

\begin{align}
(\ln \vert \psi'(x)\vert)'&=-\frac{1}{2x^{\frac{3}{2}}\sqrt{1-x}}\Big(\frac{3\sqrt{x}(1-x)-x^{\frac{3}{2}}}{\sqrt{1-x}}\Big)\nonumber\\
&=\frac{4x-3}{2x(1-x)}.\nonumber
\end{align}

\noindent This implies  $\displaystyle\frac{d}{dx}(x( \ln \vert \psi'(x)\vert)')= \frac{d}{dx}\Big(x\frac{4x-3}{2(1-x)x}\Big)=\frac{1}{2(1-x)^{2}}>0.$}
\end{enumerate}}
\item{Denote $\displaystyle h(\beta, q_j,q_n) :=\frac{\psi'(\beta q_j)}{\psi'(\beta q_n)}.$ Then $\displaystyle \frac{\partial}{\partial\beta}\Big[ \frac{h_{\beta}(\beta, q_j,q_n)}{h_{\beta}(\beta, q_l,q_n)}\Big]\neq 0$, where $ q_j <q_n <q_l.$ \\

\begin{enumerate}[(i)]
\item{
\noindent Suppose $\phi(\cdot)=\tanh(\cdot).$ Then $\displaystyle \psi'(\beta q_j)=-\frac{1}{2\beta q_j \sqrt{1-\beta q_j}}.$ This implies $\displaystyle h(\beta, q_j,q_n)=\frac{q_n}{q_j}\sqrt{\frac{1-\beta q_n}{1-\beta q_j}}.$\\
\noindent  Next, we will compute the quotient $\displaystyle \frac{h_{\beta}(\beta, q_j,q_n)}{h_{\beta}(\beta ,q_l,q_n)}.$ \\

\noindent It is easy to see that $\displaystyle h_{\beta}(\beta, q_j,q_n)=\frac{q_n}{2q_{j}}\cdot \frac{q_j -q_n}{\sqrt{1-\beta q_n}(1-\beta q_j)^{3/2}}.$ Therefore, $\displaystyle \frac{h_{\beta}(\beta, q_j,q_n)}{h_{\beta}(\beta, q_l,q_n)}=\frac{q_l(q_j -q_n)}{q_j(q_l-q_n)}\cdot \frac{(1-\beta q_l)^{3/2}}{(1-\beta q_j)^{3/2}}.$ It can be easily checked that $ \displaystyle \frac{\partial}{\partial \beta}\Big(\frac{h_{\beta}(\beta, q_j,q_n)}{h_{\beta}(\beta, q_l,q_n)}\Big)\neq 0$ for all $(\beta,q_j,q_n,q_l)$ where $q_j<q_n<q_l.$}
\item{
\noindent Now consider $\phi(\cdot)=\arctan(\cdot).$ Then $\displaystyle \psi'(\beta q)= -\frac{1}{2(\beta q)^{\frac{3}{2}}}\frac{1}{\sqrt{1-\beta q}}.$ \\
\noindent Using the expression for $\psi'(\beta q),$ we obtain  $\displaystyle h(\beta, q_j,q_n)=\Big(\frac{q_n}{q_j}\Big)^{\frac{3}{2}}\sqrt{\frac{1-\beta q_n}{1-\beta q_j}},$ and $\displaystyle h(\beta, q_l,q_n)=\Big(\frac{q_n}{q_l}\Big)^{\frac{3}{2}}\sqrt{\frac{1-\beta q_n}{1-\beta q_l}}.$ These are identical to case of $\phi(\cdot)=\tanh(\cdot).$ Remaining details can be checked in similar manner.}

\end{enumerate}}
\item{For all $p>q,$ we have $\frac{d}{d\beta}\Big(\frac{\psi(\beta p)}{\psi(\beta q)}\Big)<0.$
\begin{enumerate}[(i)]
\item{ Consider $\phi(\cdot)=\tanh(\cdot).$\\
\noindent  Then $\displaystyle \frac{\psi(\beta p)}{\psi(\beta q)}=\frac{\arctanh(\sqrt{1-\beta p})}{\arctanh(\sqrt{1-\beta q})}.$ Using definition of $\arctanh(\cdot),$ we get

$\displaystyle \frac{\psi(\beta p)}{\psi(\beta q)}=\frac{\ln\Big(\frac{1+\sqrt{1-\beta p}}{1-\sqrt{1-\beta p}}\Big)}{\ln\Big(\frac{1+\sqrt{1-\beta q}}{1-\sqrt{1-\beta q}}\Big)}.$ Now we will compute $\displaystyle \frac{d}{d\beta}\Big( \frac{\psi(\beta p)}{\psi(\beta q)}\Big).$ \\

\noindent It is easy to see that $\displaystyle \frac{d}{d\beta}\Big(\ln\Big(\frac{1+\sqrt{1-\beta p}}{1-\sqrt{1-\beta p}}\Big)\Big)=-\frac{1}{\beta\sqrt{1-\beta p}}.$\\ 
\noindent Therefore 
$\displaystyle \frac{d}{d\beta}\Big( \frac{\ln\Big(\frac{1+\sqrt{1-\beta p}}{1-\sqrt{1-\beta p}}\Big)}{\ln\Big(\frac{1+\sqrt{1-\beta q}}{1-\sqrt{1-\beta q}}\Big)}\Big)=\frac{-\frac{1}{\beta\sqrt{1-\beta p}}\ln\Big(\frac{1+\sqrt{1-\beta q}}{1-\sqrt{1-\beta q}}\Big) +\frac{1}{\beta\sqrt{1-\beta q}}\ln\Big(\frac{1+\sqrt{1-\beta p}}{1-\sqrt{1-\beta p}}\Big)}{\ln\Big(\frac{1+\sqrt{1-\beta q}}{1-\sqrt{1-\beta q}}\Big)^{2}}.$\\

\noindent Recall that $\displaystyle \arctanh(x)=\ln \Big(\frac{1+\sqrt{1-x}}{1-\sqrt{1-x}}\Big).$ Since $p>q,$ we get $\frac{1}{\sqrt{1-\beta p}}>\frac{1}{\sqrt{1-\beta q}},$ and\\
 $\arctanh (\sqrt{1-\beta p})< \arctanh(\sqrt{1-\beta q}).$ This implies that 
 $$\frac{1}{\beta\sqrt{1-\beta q}}\arctanh (\sqrt{1-\beta p}) -\frac{1}{\beta\sqrt{1-\beta p}}\arctanh (\sqrt{1-\beta q})<0.$$ 

\noindent Hence $\frac{d}{d\beta}\Big(\frac{\psi(\beta p)}{\psi(\beta q)}<0,$ where $p>q.$}
\item{ Now, consider $\phi(\cdot)=\arctan(\cdot).$\\

\noindent Then $\frac{\psi(\beta p)}{\psi(\beta q)}=\sqrt{\frac{q}{p}}\sqrt{\frac{1-\beta p}{1-\beta q}}.$ This implies that 
\begin{align}
\frac{d}{d\beta}\Big(\frac{\psi(\beta p)}{\psi(\beta q)}\Big)&= \sqrt{\frac{q}{p}}\frac{1}{2}\sqrt{\frac{1-\beta q}{1-\beta p}}\Big(\frac{-p(1-\beta q)-(-q)(1-\beta p)}{(1-\beta q)^{2}}\Big)\nonumber\\
&=\frac{1}{2}\cdot \sqrt{\frac{q}{p}}\sqrt{\frac{1-\beta q}{1-\beta p}}\Big(\frac{q-p}{(1-\beta q)^{2}}\Big)
<0\; (\text{since\; }p>q).\nonumber
\end{align}}

\end{enumerate}}

\item{For all $x>0,$ function $\psi(\cdot)$ satisfies the following property;
$$\frac{d}{dx}\Big(x\frac{d}{dx}(\frac{\psi(x)}{x \psi'(x)})\Big)\geq 0$$
\begin{enumerate}[(i)]
\item{ Consider $\phi(\cdot)=\tanh(\cdot).$ For this case $\psi(x)=\arctanh (\sqrt{1-x}),$ and $\psi'(x)=-\frac{1}{2x\sqrt{1-x}}.$ \\

We obtain
\begin{align}
&\frac{d}{dx}\Big(\frac{\psi(x)}{x \psi'(x)}\Big)= \frac{d}{dx}\Big(\frac{\arctanh (\sqrt{1-x)}}{x\cdot \frac{-1}{2x\sqrt{1-x}}}\Big)\nonumber\\
&=2\Big(\frac{\arctanh (\sqrt{1-x})}{2\sqrt{1-x}}+\frac{1}{2x}\Big)=\frac{\arctanh (\sqrt{1-x})}{\sqrt{1-x}}+\frac{1}{x}.\nonumber
\end{align}

 \noindent Therefore 
\begin{align}
&\frac{d}{dx}\Big(x\frac{d}{dx}(\frac{\psi(x)}{x \psi'(x)})\Big)=\frac{d}{dx}\Big(\frac{x \; \arctanh(\sqrt{1-x)}}{\sqrt{1-x}}+1\Big) \nonumber\\
&=\frac{1}{2(1-x)}\Big(\sqrt{1-x}\arctanh(\sqrt{1-x})+
\frac{\arctanh(\sqrt{1-x})}{\sqrt{1-x}}-1\Big).\nonumber
\end{align}

\noindent We know that $\frac{d}{dx}(\arctanh(x))>1$ for all $x$ such that $0<x<1$,  which, in turn  implies that $\displaystyle \frac{\arctanh(\sqrt{1-x})}{\sqrt{1-x}}>1.$ Hence $\frac{d}{dx}\Big(x\frac{d}{dx}(\frac{\psi(x)}{x \psi'(x)})\Big)>0.$}
\item{ Now suppose  $\phi(\cdot)=\arctan(\cdot).$ Here  $\psi(x)=\sqrt{\frac{1}{x}-1}.$ Moreover $\psi'(x)=-\frac{1}{2x^{\frac{3}{2}}\sqrt{1-x}}.$\\

\noindent We obtain
\begin{align}
&\frac{d}{dx}\Big(\frac{\psi(x)}{x\psi'(x)}\Big)=-\frac{d}{dx}\Big(\frac{\sqrt{1-x}}{x \sqrt{x}}\cdot 2\sqrt{1-x} x^{\frac{3}{2}}\Big)\nonumber\\
&=2.\nonumber
\end{align}

\noindent Now, it is easy to see that $\frac{d}{dx}\Big(x\frac{d}{dx}(\frac{\psi(x)}{x \psi'(x)})\Big)=2>0.$}

\end{enumerate}}
\end{enumerate}
\section{Conclusion}
\noindent We have studied existence of points of local maxima for the function $f(x)=\sum_{i=1}^{n}c_i \phi(x_i),$ over a hyperplane.  We have found conditions, imposed on the function $\phi(\cdot)$, which guarantee existence of at most one point of local maximum for the function on the hyperplane. Those conditions are satisfied by wide range of neuron transfer functions. The next step involves computing the points of local maxima for the nonlinear function over planes of lower dimension.

\end{document}